\documentclass[12pt,a4paper]{article}
\usepackage{amssymb}

\usepackage{graphicx}
\usepackage[tbtags,reqno]{amsmath}

\input{tcilatex}
\begin{document}

\title{The homology theory of Koszul-Vinberg algebroids and Poisson manifolds II}

\author{M. NGUIFFO BOYOM \\
D\'epartement de Math\'ematiques, Universit\'e Montpellier II \\
boyom@Math.univ-montp2.fr}


\maketitle

\begin{abstract}
We deal with smooth real manifolds as well as complex analytic manifolds as well. It
is well known that the concept of star product is powerful enough to produce all
Poisson structures on real manifolds. According to [BdM] it is not known whether
holomorphic star products exist on complex analytic manifolds.

The main purpose of this paper is to show that the concept of homology of
Koszul-Vinberg algebroids on smooth (resp. complex analytic) manifolds is an
effective tool to produce smooth (resp complex analytic) Poisson structures on
smooth (resp. complex analytic) manifolds. We also study some invariants of contact
structures which arise from the associated Koszul-Vinberg algebroids.
\end{abstract}

\textbf{Introduction. }It is well known that Lie algebroids are related to special
Poisson manifolds. On the other hand Poisson tensors appear also as infinitesimal
deformations of the associative algebras of smooth functions on real manifolds
[Kontsevich...], [DWL], this being highlighted by the concept of star product.

According to [BdM] the real arguments that are used to produce smooth star products
don't work in the case of complex analytic functions on holomorphic manifolds.

In [NGB$_{3}$] we have developped a homology theory of Koszul-algebras and their
modules. That homology provides a framework for many differential geometry problems
such as deformations of hyperbolic affine manifolds. The purpose of this paper is to
use the homology theory of Koszul-Vinberg algebroids to produce all Poisson
structures on smooth (resp. complex analytic) manifolds.

To every smooth (resp. complex analytic) manifold is assigned a canonical
Koszul-Vinberg algebroid $\mathcal{V}\left( M\right) $ whose base manifold is the
cotangent bundle (resp. the holomorphic cotangent manifold) $T^{\ast }M.$ We focus
on real (complex) homology of $\mathcal{V}\left( M\right) ,$
say $H_{\ast }\left( \mathcal{V}\left( M\right) ,\mathbb{R}\right) $ (resp. $%
H_{\ast }\left( \mathcal{V}\left( M\right) ,\mathbb{C}\right) ).$ We show that there
is a one to one correspondence between the set of smooth \ (resp. holomorphic)
Poisson manifold structures on M, say $\mathcal{P}\left(
M\right) ,$ and a subset of $H_{2}\left( \mathcal{V}\left( M\right) ,\mathbb{%
R}\right) $ (resp. $H_{2}\left( \mathcal{V}\left( M\right) ,\mathbb{C}%
\right) ).$ We also point out a one to one correspondence between the vector
space of smooth (resp. holomorphic) vector fields on $M$ and a subspace of $%
H_{1}\left( \mathcal{V}\left( M\right) ,\mathcal{V}\left( M\right) \right) .$

\bigskip
\bigskip

\qquad \textbf{1 - Backgroud materials}

\bigskip

\qquad Throughout this section we denote by $\mathbb{F}$ a commutative field
of characteristic zero. All vector spaces that are considered have the same
base field, $\mathbb{F}.$

\bigskip

\textbf{1.1. Koszul-Vinberg algebras and their modules.}

\qquad Given an algebra $\mathcal{A}$ whose multiplication map is denoted by

\begin{equation*}
ab,\quad a\in \mathcal{A},b\in \mathcal{A}.
\end{equation*}

the associator map of $\mathcal{A}$ is the 3-multilinear map, $\left(
,,\right) $ defined by

\begin{equation}
\left( a,b,c\right) =a\left( bc\right) -\left( ab\right) c,\quad \forall
a,b,c\in \mathcal{A}  \tag{1}
\end{equation}

Let $W$ be a vector space endoved with two $\mathcal{A}$-actions

\begin{eqnarray*}
\mathcal{A\times }W &\rightarrow &W \\
\left( a,w\right) &\rightarrow &aw
\end{eqnarray*}

and

\begin{eqnarray*}
W\times \mathcal{A} &\rightarrow &W \\
\left( w,a\right) &\rightarrow &wa.
\end{eqnarray*}

Given elements $a,b\in \mathcal{A}$ and $w\in W$ we set

\begin{eqnarray}
\left( a,b,w\right) &=&a\left( bw\right) -\left( ab\right) w  \TCItag{2} \\
\left( a,w,b\right) &=&a\left( wb\right) -\left( aw\right) b,  \notag \\
\left( w,a,b\right) &=&w\left( ab\right) -\left( wa\right) b.  \notag
\end{eqnarray}

\textbf{D\'{e}finition. }\textit{An algebra }$\mathcal{A}$ \textit{is called
a Koszul-Vinberg algebra if its associator map is symmetric w.r.t. the first
two arguments.}

\bigskip

Given a Koszul-Vinberg algebra $\mathcal{A},$ a vector space $W$ with two $%
\mathcal{A}$-actions $\mathcal{A\times }W\rightarrow W$ and $W\times
\mathcal{A\rightarrow }W$ we have :

\bigskip

\textbf{D\'{e}finition. }\textit{A vector space }$W$\textit{\ with two
actions as above is called a Koszul-Vinberg module of }$\mathcal{A}$\textit{%
\ if it\ satisfies the following identities :}

\begin{eqnarray*}
\left( a,b,w\right) &=&\left( b,a,w\right) , \\
\left( a,w,b\right) &=&\left( w,a,b\right) . \\
\forall a,\forall b &\in &\mathcal{A}\text{ and }\forall w\in W
\end{eqnarray*}

Here are examples of Koszul-Vinberg algebras.

$\left( e_{1}\right) .$\quad Every associative algebra is Koszul-Vinberg
algebra.

$\left( e_{2}\right) .$\quad Let us consider $\mathbb{F}^{4}=\left\{
X=\left( x,y,z,t\right) ,\text{ Given }x,y,z,t\in \mathbb{F}\right\} $

\begin{eqnarray*}
X &=&\left( x,y,z,t\right) , \\
X^{\prime } &=&\left( x^{\prime },y^{\prime },z^{\prime },t^{\prime }\right)
\end{eqnarray*}

we define the multiplication

\begin{equation*}
XX^{\prime }=\left( \left( y+t\right) z^{\prime }-tt^{\prime },zz^{\prime
}-t\left( x^{\prime }+t^{\prime }\right) ,\left( y-t\right) t^{\prime
},0\right) .
\end{equation*}

It is easy to verify that $\forall X,X^{\prime },X"\in \mathbb{F}^{4}$ we
have

\begin{equation*}
\left( X,X^{\prime },X"\right) =\left( X^{\prime },X,X"\right) .
\end{equation*}

$\left( e_{3}\right) .\quad $Let $\mathbb{F}$ be the field $\mathbb{R}$ of
real numbers. Let $L$\ be the vector space of formal vector fields at the
origin of $\mathbb{R}^{m}.$ Thus $\forall X\in L,$

\begin{eqnarray*}
X &=&\sum\limits_{i=0}^{n}X_{i}\frac{\partial }{\partial x_{i}}, \\
X_{i} &\in &\mathbb{R}\left[ \left[ x_{1},...,x_{m}\right] \right] .
\end{eqnarray*}

Given two elements of $L,$ say $X$\ and $X^{\prime }$\ we set

\begin{equation*}
XX^{\prime }=\sum_{i=1}^{m}\left( \sum_{j=1}^{m}X_{i}\frac{\partial
X_{j}^{^{\prime }}}{\partial x_{i}}\right) \frac{\partial }{\partial x_{i}}
\end{equation*}

It is easy to check that L. with the above multiplication is a
Koszul-Vinberg algebra

\begin{equation*}
\left( X,X^{\prime },X"\right) =\left( X^{\prime },X,X"\right) .
\end{equation*}

Let $V,W$ be two Koszul-Vinberg modules of a Koszul-Vinberg algebra $%
\mathcal{A}.$ The (algebraic) tensor product $V\otimes W$ is also a
Koszul-Vinberg module of $\mathcal{A}$ under the following actions

\begin{eqnarray*}
a\left( v\otimes w\right) &=&av\otimes w+v\otimes aw \\
\left( v\otimes w\right) a &=&v\otimes wa
\end{eqnarray*}

In particular the vector spaces $\ T^{r}V=\bigotimes\limits^{r}V,r\in
\mathbb{N}^{\ast },$ \ are Koszul-Vinberg modules of $\ \mathcal{A}.$ \ We
will equip the vector space \ $\hom _{\mathbb{F}}\left( T^{r}W,T^{s}V\right)
$ \ \ with the following $\ \ \mathcal{A}$-actions :

For given $\theta \in $Hom$_{\mathbb{F}}\left( T^{r}W,T^{s}V\right) ,a\in
\mathcal{A},$ we set

\begin{eqnarray}
\left( a\theta \right) \left( \xi \right) &=&a\left( \theta \left( \xi
\right) \right) -\theta \left( a\xi \right) ,  \TCItag{3} \\
\left( \theta a\right) \left( \xi \right) &=&\left( \theta \left( \xi
\right) \right) a\quad \forall \xi \in T^{r}W.  \notag
\end{eqnarray}

Under (3) Hom$_{\mathbb{F}}\left( T^{r}W,T^{s}V\right) $ is a Koszul-Vinberg
module of $\mathcal{A}$.

Given a Koszul-Vinberg module $W$ of a Koszul-Vinberg algebra $\mathcal{A},$
we set

\begin{equation*}
J\left( W\right) =\left\{ w\in W/\left( a,b,w\right) =0\quad \forall
a,\forall b\in \mathcal{A}\right\}
\end{equation*}

In particular regarding $\mathcal{A}$ as Koszul-Vinberg module of itself we
set

\begin{equation*}
J\left( A\right) =\left\{ c\in \mathcal{A}/\left( a,b,c\right) =0\quad
\forall a,\forall b\in \mathcal{A}\right\} .
\end{equation*}

\bigskip

\textbf{1.2 - The chain complex }$C\left( \mathcal{A},W\right) .$

\qquad To every Koszul-Vinberg module of $\mathcal{A},$ say $W,$ we assign
the $\mathbb{Z}$-graded vector space

\begin{equation*}
C\left( \mathcal{A},W\right) =\bigoplus_{q\in \mathbb{Z}}C_{q}\left(
\mathcal{A},W\right)
\end{equation*}

where

\begin{eqnarray*}
C_{q}\left( \mathcal{A},W\right) &=&0\text{ if }q<0, \\
C_{0}\left( \mathcal{A},W\right) &=&J\left( W\right) , \\
C_{q}\left( \mathcal{A},W\right) &=&\text{Hom}_{\mathbb{F}}\left( T^{q}%
\mathcal{A},W\right) \text{ if }q>0
\end{eqnarray*}

Let $q\in \mathbb{Z},q>0.$ For every $\left( a,j\right) \in \mathcal{A\times
}\mathbb{Z}$ with $1\leqslant j\leqslant q$ we define a linear map $%
e_{j}\left( a\right) $ from $C_{q}\left( \mathcal{A},W\right) $ to $%
C_{q-1}\left( \mathcal{A},W\right) $ by putting

\begin{equation*}
e_{j}\left( a\right) \theta \left( a_{1}\ldots a_{q-1}\right) =\theta \left(
a_{1}\ldots a_{j-1},a,a_{j},\ldots a_{q-1}.\right)
\end{equation*}

We are going now to equip the $\mathbb{Z}$-graded vector space $C\left(
\mathcal{A},W\right) $ with a differential operator of degree +1. Indeed we\
define

\begin{equation*}
\delta _{q}:C_{q}\left( \mathcal{A},W\right) \rightarrow C_{q+1}\left(
\mathcal{A},W\right)
\end{equation*}

by

\begin{equation}
\delta _{0}\left( w\right) \left( a\right) =-aw+wa\quad \forall a\in
\mathcal{A},\forall w\in J\left( W\right) .  \tag{3.1}
\end{equation}

If $q>0$

\begin{eqnarray}
\left( \delta _{q}\theta \right) \left( a_{1}\ldots a_{q+1}\right)
&=&\sum_{j\leqslant q}\left( -1\right) ^{j}\left\{ \left( a_{j}\theta
\right) \left( \ldots \hat{a}_{j}\ldots a_{q+1}\right) \right.  \TCItag{3.2}
\\
&&\left. +\left( e_{q}\left( a_{j}\right) \theta a_{q+1}\right) \left(
\ldots \hat{a}_{j}\ldots \hat{a}_{q+1}\right) \right\}  \notag
\end{eqnarray}

The operator $\delta =\left( \delta _{q}\right) _{q}$ satisfies the
following identity

\begin{equation*}
\delta _{q+1}\circ \delta _{q}=0
\end{equation*}

Details can be found in $\left[ \text{NGB}_{3}\right] .$ Thus we get the
chain complex

\begin{equation}
J\left( W\right) \overset{\delta _{0}}{\rightarrow }C_{1}\left( \mathcal{A}%
,W\right) \rightarrow \ldots C_{q}\left( \mathcal{A},W\right) \overset{%
\delta q}{\rightarrow }\ldots  \tag{4}
\end{equation}

The $q^{th}$ homology of $\left( 4\right) $ is denoted by $H_{q}\left(
\mathcal{A},W\right) :$

\begin{equation*}
H_{q}\left( \mathcal{A},W\right) =\ker \left( \delta _{q}\right) /\text{im}%
\left( \delta _{q-1}\right) .
\end{equation*}

It is to be noticed that the chain complex $\left( 4\right) $ is different
from that constructed by Albert Nijenhuis in [N.A]. In [NGB$_{3}$] we have
shown that the complex $\left( 4\right) $ provides an efficient framework
for the deformation theory of Koszul-Vinberg algebras and Koszul-Vinberg
modules. Thus from the viewpoint of deformation theory the complex $\left(
4\right) $ answers the fundamental question raised by M. Gerstenhaber ; \
viz ``every restricted deformation theory generates its proper cohomology
theory'', [GM]. In the litterature Koszul-Vinberg algebras are \ also called
left symmetric algelbras [NGB$_{1}$], [HJ], [PM], they are closely related
to differential geometry of bounded domains and to affinaly flat geometry,
[VEB], [MJ], [NGB$_{4}].$

\bigskip

\textbf{1.3. Extension of affine structure and }$H_{2}\left( \mathcal{A},%
\mathbb{F}\right) .$

\qquad Let $\mathcal{A}$ be a Koszul-Vinberg algebra ; it gives rise to the
Lie algebra $A_{L}$ whose bracket is defined by

\begin{equation*}
\left[ a,b\right] =ab-ba.
\end{equation*}

Let us denote by $\Pi $ the linear map from Hom$_{\mathbb{F}}\left( T^{2}%
\mathcal{A},\mathcal{A}\right) $ to Hom$_{\mathbb{F}}\left( \overset{2}{%
\Lambda }\mathcal{A},\mathcal{A}\right) $ which is defined by

\begin{equation*}
\left( \Pi \theta \right) \left( a,b\right) =\frac{1}{2}\left( \theta \left(
a,b\right) -\theta \left( b,a\right) \right) .
\end{equation*}

We consider $\mathbb{F}$ as a trivial module of $\mathcal{A},$ it becomes a
trivial module of $\mathcal{A}_{L},$ thus we will consider the Hochschild
complex

\begin{equation}
\ldots \rightarrow C^{q}\left( \mathcal{A}_{L},\mathbb{F}\right) \overset{%
\delta }{\rightarrow }C^{q+1}\left( \mathcal{A}_{L},\mathbb{F}\right)
\rightarrow ...  \tag{5}
\end{equation}

The $\mathbb{Z}$-graded vector space

\begin{equation*}
C^{\ast }\left( \mathcal{A,}\mathbb{F}\right) =\bigoplus_{q\in \mathbb{Z}%
}C^{q}\left( \mathcal{A}_{L},\mathbb{F}\right)
\end{equation*}

is defined by

\begin{eqnarray*}
C^{q}\left( \mathcal{A}_{L},\mathbb{F}\right) &=&0\text{ if }q<0, \\
C^{0}\left( \mathcal{A}_{L},\mathbb{F}\right) &=&\mathbb{F} \\
C^{q}\left( \mathcal{A}_{L},\mathbb{F}\right) &=&\text{Hom}\left( \overset{q%
}{\Lambda }\mathcal{A}_{L},\mathbb{F}\right) \text{ if }q>0
\end{eqnarray*}

Thus

\begin{eqnarray*}
\delta _{0} &=&0\text{ and} \\
\delta \theta \left( a_{1}\ldots a_{q+1}\right) &=&\sum_{i<j}\left(
-1\right) ^{i+j}\theta \left( \left[ a_{i},a_{j}\right] ,\ldots \hat{a}%
_{i}\ldots \hat{a}_{j}\ldots \right)
\end{eqnarray*}
Now let us regard $\Pi $ as map from $C_{2}\left( \mathcal{A},\mathbb{F}%
\right) $ to $C^{2}\left( \mathcal{A}_{L},\mathbb{F}\right) .$ Then one
easily verifies that

\begin{equation*}
\oint \left( \delta _{2}\theta \right) \left( a,b,c\right) =\left( 2\delta
\Pi \theta \right) \left( a,b,c\right)
\end{equation*}

where $\oint $ stands for cyclic sum w.r.t. a,b,c. Thus $\Pi $ induces a
linear map

\begin{equation*}
\Pi :H_{2}\left( \mathcal{A},\mathbb{F}\right) \rightarrow H^{2}\left( A_{L},%
\mathbb{F}\right) .
\end{equation*}

By setting

\begin{equation*}
\mathcal{O}\left( \mathcal{A},\mathcal{A}_{L}\right) =H^{2}\left( \mathcal{A}%
_{L},\mathbb{F}\right) /\text{im}\left( \Pi \right)
\end{equation*}

we get the exact sequence of vector spaces

\begin{equation}
0\rightarrow \text{Ker}\left( \Pi \right) \rightarrow H_{2}\left( \mathcal{A}%
,\mathbb{F}\right) \rightarrow H^{2}\left( \mathcal{A}_{L},\mathbb{F}\right)
\overset{\sigma }{\rightarrow }\mathcal{O}\left( A,A_{L}\right) \rightarrow 0
\tag{6}
\end{equation}

The quotient space $\mathcal{O}\left( \mathcal{A}_{L}\right) $ has a
geometric meaning. Indeed let $\omega $ be cocycle in $C^{2}\left( \mathcal{A%
}_{L},\mathbb{F}\right) .$ It defines a Lie algebra structure in $\mathbb{%
F\oplus }\mathcal{A}_{L}$ whose bracket is given by

\begin{equation*}
\left[ \left( \lambda ,a\right) ,\left( \mu ,b\right) \right] =\left( \omega
\left( a,b\right) ,\left[ a,b\right] \right) .
\end{equation*}

Since the bracket of $\mathcal{A}_{L},$ has the form

\begin{equation*}
\left[ a,b\right] =ab-ba
\end{equation*}

where $\left( a,b\right) \rightarrow ab$ is the multiplication in the
Koszul-Vinberg algebra, the question arises to know whether $\mathbb{F}%
\bigoplus\limits_{\omega }\mathcal{A}_{L}$ admits a Koszul-Vinberg
multiplication

\begin{equation*}
\left( \left( \lambda ,a\right) ,\left( \mu ,b\right) \right) \rightarrow
\left( \lambda ,a\right) \cdot \left( \mu ,b\right)
\end{equation*}

such that

\begin{equation*}
\left( \lambda ,a\right) \cdot \left( \mu ,b\right) -\left( \mu ,b\right)
\cdot \left( \lambda ,a\right) =\left( \omega \left( a,b\right) ,\left[ a,b%
\right] \right)
\end{equation*}

If such a multiplication exists then by setting

\begin{equation*}
\left( 0,a\right) \cdot \left( 0,b\right) =\left( \varphi \left( a,b\right)
,ab\right)
\end{equation*}

we see that $\delta _{2}\varphi =0$ and $2\Pi \varphi \left( a,b\right)
=\omega \left( a,b\right) .$ Thus

\begin{equation*}
\sigma \left( \left[ \omega \right] \right) =0
\end{equation*}

where $\left[ \omega \right] $ is the class of $\omega $ in $H^{2}\left(
\mathcal{A}_{L},\mathbb{F}\right) .$

We conclude that given a 2-cocycle $\omega $ in $C^{2}\left( \mathcal{A}_{L},%
\mathbb{F}\right) ,$ the Lie algebra $\mathbb{F}\bigoplus\limits_{\omega }%
\mathcal{A}_{L}$ admits a Koszul-Vinberg algebra structure.

\begin{equation*}
\left( \lambda ,a\right) \underset{\varphi }{\cdot }\left( \mu ,b\right)
=\left( \varphi \left( a,b\right) ,ab\right)
\end{equation*}

such that

\begin{equation*}
\left( \lambda ,a\right) \underset{\varphi }{\cdot }\left( \mu ,b\right)
-\left( \mu ,b\right) \underset{\varphi }{\cdot }\left( \lambda ,a\right)
=\left( \omega \left( a,b\right) ,\left[ a,b\right] \right)
\end{equation*}

if and only if $\sigma \left( \left[ \omega \right] \right) =0.$ This
obstruction meaning of $\mathcal{O}\left( \mathcal{A}_{L},\mathcal{A}\right)
$ is a way to understand the example of nilpotent Lie groups without left
invariant affine structure, $\left[ BY\right] .$ To be more explicit let us
consider a finite dimensional real Koszul-Vinberg algebra $\mathcal{A}.$ Let
$G$ be the connected and simply connected Lie group whose Lie algebra is $%
\mathcal{A}_{L}.$ Thus the multiplication map in $\mathcal{A}$ gives rise to
a locally flat linear connection on $G,$ say $D$ which is invariant under
the left translations by elements of $G.$ Given $X,X^{\prime }\in \mathcal{A}%
_{L}$ we have

\begin{equation*}
D_{X}X^{\prime }=XX^{\prime }.
\end{equation*}

Let $\omega $ be a cocycle in $C^{2}\left( \mathcal{A}_{L},\mathbb{R}\right)
,$ then $\omega $ defines the following bracket in $\frak{G}_{\omega }=%
\mathbb{R\oplus }\mathcal{A}_{L}.$

\begin{equation}
\left[ \left( \lambda ,X\right) ,\left( \lambda ^{\prime },X^{\prime
}\right) \right] =\left( \omega \left( X,X^{\prime }\right) ,\left[
X,X^{\prime }\right] \right) .  \tag{7}
\end{equation}

So, we have the exact sequence of Lie algebra

\begin{equation}
0\rightarrow \mathbb{R\rightarrow }\frak{G}_{\omega }\rightarrow \mathcal{A}%
_{L}\frak{\rightarrow 0}  \tag{8}
\end{equation}

\qquad Let $G_{\omega }$ be the connected and simply connected real Lie
group whose Lie algebra is $\frak{G}_{\omega },$ then $\left( 8\right) $ is
equivalent to the exact sequence of real Lie groups

\begin{equation}
1\rightarrow \mathbb{R}^{>0}\rightarrow G_{\omega }\overset{p}{\rightarrow }%
G\rightarrow 1.  \tag{8'}
\end{equation}

Now the question arises to know whether $G_{\omega }$ does admit a left
invariant locally linear connection $D_{\omega }$ such that the projection $%
p $ is an affine homomorphism from $\left( G_{\omega },D_{\omega }\right) $
to $\left( G,D\right) .$ The obstruction to the existence of $D_{\omega }$is
the image of the class $\left[ \omega \right] \in H^{2}\left( \mathcal{A}%
_{L},\mathbb{R}\right) $ under the linear map

\begin{equation*}
\sigma :H^{2}\left( \mathcal{A}_{L},\mathbb{R}\right) \rightarrow \mathcal{O}%
\left[ \mathcal{A}_{L},\mathcal{A}\right]
\end{equation*}

\bigskip

\textbf{2 - Scalar homology of Koszul-Vinberg algebroids.}

\bigskip

\qquad The aim of this section is to introduce the concept of real (resp.
complex) homology of smooth (resp. complex analytic) Koszul-Vinberg
algebroids according to $\mathbb{F=R}$ or $\mathbb{F=C.}$

\bigskip

\textbf{2.1. Koszul-Vinberg algebroids.}

\qquad According to $\mathbb{F=R}$ or $\mathbb{F=C},$ manifold $M$ will mean
for smooth or complex analytic manifold. Vector bundles over $M$ and
sections of vector bundles are smooth or holomorphic according to $\mathbb{%
F=R}$ or $\mathbb{F=C}.$ In particular if $M$ is holomorphic manifold then $%
\Gamma \left( TM\right) $ stands for holomorphic vector fields. So given a
vector bundle over $M,$ say $E,\Gamma \left( E\right) $ is the vector space
of sections of $E\rightarrow M$ ; $\Gamma \left( \overset{\sim }{\mathbb{F}}%
\right) $ stands for the associative algebra of smooth (resp. holomorphic)
functions on $M$ according to $\mathbb{F=R}$ or $\mathbb{F=C}.$

\bigskip

\textbf{Definition}\textit{. A Koszul-Vinberg algebroid over a manifold }$M$%
\textit{\ is a couple }$\left( E,a_{n}\right) $ where $E$ is a vector bundle
over $M$ \textit{and }$a_{n}$\textit{\ is a }$\mathbb{F}$\textit{-linear map
(anchor) from }$\Gamma \left( E\right) $\textit{\ to }$\Gamma \left(
TM\right) $\textit{\ satisfying the following conditions :}

$\left( r_{1}\right) \quad \Gamma \left( E\right) $\textit{\ is a
Koszul-Vinberg algebra (,its multiplication will be denoted by }$ss^{\prime
}\quad \forall s,$ $s^{\prime }\in \Gamma \left( E\right) )$ ;

$\left( r_{2}\right) \quad $\textit{given }$f\in \Gamma \left( \overset{\sim
}{\mathbb{F}}\right) ,s\in \Gamma \left( E\right) ,s^{\prime }\in \Gamma
\left( E\right) $ we have

\begin{eqnarray*}
\left( fs\right) s^{\prime } &=&f\left( ss\right) \ ; \\
s\left( fs^{\prime }\right) &=&f\left( ss^{\prime }\right) +\left(
a_{n}\left( s\right) \cdot f\right) s^{\prime }.
\end{eqnarray*}

\textbf{Examples of Koszul-Vinberg algebroids.}

$\left( e_{1}\right) \quad $Let $\left( M,D\right) $ be a locally flat
manifold. Then $\Gamma \left( TM\right) $ is a Koszul-Vinberg algebra whose
multiplication map is defined by

\begin{equation*}
\left( s,s^{\prime }\right) \rightarrow ss^{\prime }=D_{s}s^{\prime }
\end{equation*}

Thus the couple ($TM,$ identity) is Koszul-Vinberg algebroid (when $\mathbb{%
F=C},TM$ stands for $T^{1,0}\left( M\right) ).$

$\left( e_{2}\right) \quad $Let $\mathcal{L}$ be Lagrangian foliation in a
symplectic manifold $\left( M,\omega \right) .$ Let $E\left( \mathcal{L}%
\right) \subset TM$ be the tangent bundle of $\mathcal{L}$. We consider two
sections $s$ and $s^{\prime }$ of $E\left( \mathcal{L}\right) \rightarrow M$
and we define the section $ss^{\prime }\in \Gamma \left( E\left( \mathcal{L}%
\right) \right) $ by the formula

\begin{equation*}
i\left( ss^{\prime }\right) \omega =L_{s}i\left( s^{\prime }\right) \omega
\end{equation*}

where $L_{s}$ (resp. $i\left( s^{\prime }\right) )$ stands for the Lie
derivation (resp. inner product). It is easily verified that

\begin{equation*}
\left( s,s^{\prime },s"\right) =\left( s^{\prime },s,s"\right) \quad \forall
s,s^{\prime },s"\in \Gamma \left( E\left( \mathcal{L}\right) \right) .
\end{equation*}

The couple ($E\left( \mathcal{L}\right) ,$ inclusion map) is a
Koszul-Vinberg algebroid.

$\left( e_{3}\right) \quad $Particular cases of $\left( e_{2}\right) $ are
determined by completely integrable hamiltonian systems in (compact)
symplectic manifolds. Given such a system $\Sigma $ in $\left( M,\omega
\right) $ the orbits of $\Sigma $ are orbits of a locally free action of $%
\mathbb{F}^{m}$ on $M,$ $2m=\dim _{\mathbb{F}}M.$ Thus these orbits are
affine manifolds. So $\Sigma $ gives rise to a Koszul-Vinberg algebroid.

$\left( e_{4}\right) $\quad Every non singular section $X\in \Gamma \left(
TM\right) $ defines a Koszul-Vinberg algebroid

\begin{equation*}
E=\bigcup_{x\in M}\mathbb{F}X_{x}.
\end{equation*}

Sections of $E$ have the form $s=fX$ with $f\in \Gamma \left( \mathbb{F}%
\right) .$ Thus given $s=fX$ and $s^{\prime }=f^{\prime }X$ we set

\begin{equation*}
ss^{\prime }=f\left( Xf^{\prime }\right) .
\end{equation*}

Given a Koszul-Vinberg algebroid $\left( E,a_{n}\right) $ we consider the
Whitney sum $\varepsilon =E\oplus \overset{\sim }{\mathbb{F}}$ where $%
\overset{\sim }{\mathbb{F}}$ stands for the trivial bundle

\begin{equation*}
M\times \mathbb{F\rightarrow M.}
\end{equation*}

Sections of $\varepsilon $ have the form

\begin{equation*}
\xi =\left( s,f\right) \in \Gamma \left( E\right) \times \Gamma \left(
\overset{\sim }{\mathbb{F}}\right)
\end{equation*}

Let $\xi =\left( s,f\right) ,$ $\xi ^{\prime }=\left( s^{\prime },f^{\prime
}\right) $ be elements of $\Gamma \left( \varepsilon \right) ,$ we define
the section $\xi \xi ^{\prime }\in \Gamma \left( \varepsilon \right) $ by

\begin{equation}
\xi \xi ^{\prime }=\left( ss^{\prime },ff^{\prime }+a_{n}\left( s\right)
f^{\prime }\right) .  \tag{9}
\end{equation}

By direct calculation one sees that if $\xi ,\xi ^{\prime },\xi "$ are
elements of $\Gamma \left( \varepsilon \right) $ then

\begin{equation*}
\left( \xi ,\xi ^{\prime },\xi "\right) =\left( \xi ^{\prime },\xi ,\xi
"\right) .
\end{equation*}

Now let us set

\begin{equation*}
a_{\varepsilon }\left( s,f\right) =a_{n}\left( s\right) .
\end{equation*}

Then $\left( \varepsilon ,a_{\varepsilon }\right) $ is a Koszul-Vinberg
algebroid. We shall denote by $\mathcal{G}$ the vector space $\Gamma \left(
\varepsilon \right) $ endoved with the multiplication $\left( 9\right) .$

The vector space $W=\Gamma \left( \overset{\sim }{\mathbb{F}}\right) $ is a
two-sided ideal of the Koszul-Vinberg algebra $\mathcal{G}.$ Therefore we
shall regard $W$ as Koszul-Vinberg module of $\mathcal{G}\frak{.}$

Henceforth we plan dealing with the chain complexes of $E.$

\bigskip

\textbf{2.2. Homology of }$E$

\qquad Let us consider the following chain complexes :

\begin{eqnarray*}
C\left( \mathcal{G},W\right) &=&\left( C_{q}\left( \mathcal{G}\frak{,}%
W\right) ,\delta _{q}\right) ,\quad q\in \mathbb{Z}\text{,} \\
C\left( \mathcal{G}\frak{,\mathcal{G}}\right) &=&\left( C_{q}\left( \mathcal{%
G}\frak{,\mathcal{G}}\right) ,\delta _{q}\right) ,\quad q\in \mathbb{Z.}
\end{eqnarray*}

The boundary operators $\delta _{q}$ are defined by $\left( 3.1\right) $ and
$\left( 3.2\right) .$ The $q^{th}$ homology of $C\left[ \mathcal{G}\frak{,}W%
\right] $ is denoted by $H_{q}\left( M,\mathbb{F}\right) $ and that of $%
C\left( \mathcal{G},\mathcal{G}\right) $ by $H_{q}\left( \mathcal{G}\right)
. $

\textbf{Definition. }\textit{(i) Given a Koszul-Vinberg algebroid }$%
E\rightarrow M,$ \textit{the} $q^{th}$ \textit{scalar homology of }$E,$
\textit{denoted by} $H_{q}\left( E,\mathbb{F}\right) $ \textit{is the vector
space} $H_{q}\left( g,\mathbb{F}\right) .$\textit{\ (ii) The }$q^{th}$
\textit{homology of} $E,$ \textit{say} $H_{q}\left( E\right) $ \textit{is
the vector space }$H_{q}\left( \mathcal{G},\mathcal{G}\right) .$

\bigskip

We just observed that $W$ is two-sided ideal of $\mathcal{G},$ so we have
the exact sequence of Koszul-Vinberg algebras

\begin{equation}
0\rightarrow W\rightarrow \mathcal{G}\rightarrow \frak{\Gamma }\left(
E\right) \rightarrow 0  \tag{10}
\end{equation}

Since $\Gamma \left( E\right) $ is a subalgebra of $\mathcal{G}$, sequence $%
\left( 10\right) $ is splittable. On the other hand we see that $\Gamma
\left( E\right) $ is a left ideal of $\mathcal{G}$ because

\begin{equation*}
\left( s,f\right) \cdot \left( s^{\prime },0\right) =\left( ss^{\prime
},0\right) .
\end{equation*}

\bigskip

\textbf{2.3. A vanishing theorem.}

\qquad We are concerned with the homology space $H\left( E\right)
=\bigoplus\limits_{q}H_{q}\left( E\right) $ of a Koszul-Vinberg algebroid $%
\left( E,a_{n}\right) .$

For every non negative integer $q\in \mathbb{Z}$ the vector space $%
C_{q}\left( \mathcal{G}\text{,W}\right) $ is bigraded by its subspaces

\begin{equation}
C_{r,s}\left( \mathcal{G}\text{,W}\right) =Hom_{\mathbb{F}}\left( T^{r}%
\mathcal{A\otimes T}^{s}W,W\right)  \tag{11}
\end{equation}

In other words we have

\begin{equation}
C_{q}\left( \mathcal{G}\text{,W}\right)
=\bigoplus\limits_{r+s=q}C_{r,s}\left( \mathcal{G}\text{,W}\right) .
\tag{12}
\end{equation}

\qquad The bigraduation (11) does not agree with the boundary operators $%
\delta _{q}.$ However we can equip the homology space $H_{q}\left( E,\mathbb{%
F}\right) $ with the graduation defined by

\begin{equation*}
H_{r,s}\left( E,F\right) =\frac{\ker \left( \delta _{q}:C_{r,s}\left(
\mathcal{G}\text{,W}\right) \rightarrow C_{r+1,s}\left( \mathcal{G}\text{,W}%
\right) \oplus C_{r,s+1}\left( \mathcal{G}\text{,W}\right) \right) }{\delta
_{q-1}\left( C_{r-1},s\left( \mathcal{G}\text{,W}\right) +C_{r,s-1}\left(
\mathcal{G}\text{,W}\right) \right) \cap C_{r,s}\left( \mathcal{G}\text{,W}%
\right) }.
\end{equation*}

So we see that

\begin{equation*}
H_{q}\left( E,\mathbb{F}\right) =\bigoplus\limits_{r+s=q}H_{r,s}\left( E,%
\mathbb{F}\right) .
\end{equation*}

\bigskip

\textbf{Theorem I. }\textit{For every positive integer }$r$\textit{\ we have}%
$H_{r,0}\left( E,\mathbb{F}\right) =0.$

\bigskip

\textit{Proof. Let }$\theta \in C_{q}\left( \mathcal{G}\text{,W}\right) $ be
a cycle, viz $\delta _{q}\theta =0.$ We will decompose $\theta $ according
to $\left( 11\right) ,$ that is to say

\begin{equation*}
\theta =\sum_{r+s=q}\theta _{r,s}
\end{equation*}

with $\theta _{r,s}\in C_{r,s}\left( \frak{G}\text{,W}\right) .$ Since we
have

\begin{equation*}
\delta _{q}\theta _{r,s}\in C_{r+1,s}\left( \frak{G}\text{,W}\right) \oplus
C_{r,s+1}\left( \frak{G}\text{,W}\right)
\end{equation*}

we see that for $f\in W$ and $\xi _{1},\ldots ,\xi _{q}\in \mathcal{A}%
=\Gamma \left( E\right) $ we have

\begin{eqnarray*}
\delta _{1}\theta \left( f,\xi _{1\ldots }\xi _{q}\right) &=&-f\theta \left(
\xi _{1},\ldots ,\xi _{q}\right) -\sum_{1\preccurlyeq j<q}\left( -1\right)
^{j}\left( \xi _{j}\theta \right) \left( f\ldots \overset{\symbol{94}}{\xi }%
_{j}\ldots \xi _{q}\right) \\
&=&-f\theta \left( \xi _{1},\ldots ,\xi _{q}\right) -\sum_{1\preccurlyeq
j<q}\xi _{j}\theta _{q-1,1}\left( f\ldots \overset{\symbol{94}}{\xi }%
_{j}\ldots \xi _{q}\right)
\end{eqnarray*}

Let us take $f$ to be the constant function $1,$ then

\begin{equation*}
\delta _{q}\theta \left( 1,\xi _{1},\ldots ,\xi _{q}\right) =-\theta
_{q,0}\left( \xi _{1},\ldots ,\xi \right) -\delta _{q}\left( e_{1}\left(
1_{e}\right) \theta _{q-1,1}\right) \left( \xi _{1}\ldots \xi _{q}.\right)
\end{equation*}

In particular if $\delta _{q}\theta =0,$ then we obtain

\begin{equation*}
\theta _{q,0}=-\delta _{q}e_{1}\left( 1_{e}\right) \theta _{q-1,1}.
\end{equation*}

This ends the proof of Theorem I \ $\bullet $

\bigskip

\textbf{3.- Symbols of chains of superorder.}

\bigskip

\qquad This section is inspired by the concept of homology of complex of
differential forms of super order by J.L. Koszul, [KJL], (see also the
Spencer homology theory, [SS], [GH], [MB].

\bigskip

\textbf{3.1. The canonical decomposition of homogeneous chains.}

\qquad We keep all conventions of Section 2.

Consider a Koszul-Vinberg algebroid $\left( E,a_{n}\right) $ over a manifold
$M.$ We associate to $E$ the Koszul-Vinberg algebra $\mathcal{G=\Gamma }%
\left( E\oplus \overset{\sim }{\mathbb{F}}\right) $ whose multiplication is
given by $\left( 9\right) .$

\bigskip

\textbf{Definition.}\textit{Given a non negative integer }$k$\textit{\ a
homogeneous chain }$\theta \in C_{q}\left( \mathcal{G},W\right) $ \textit{is
said to be of order }$\leqslant k$ \textit{if for }$\left( \xi _{1},\ldots
,\xi _{q}\right) \in \mathcal{G}^{q}$ \textit{the value at each }$x\in M$%
\textit{\ of }$\theta \left( \xi _{1},\ldots ,\xi _{q}\right) $\textit{\
depends on the k}$^{th}$\textit{\ jet at }$x$\textit{\ of the }$\xi
_{j}^{^{\prime }s}.$

\bigskip

From now on we shall present $j_{x}^{k}\xi $ as follows :

\begin{equation*}
j_{x}^{k}\xi =\left( \xi \left( x\right) ,d_{x}^{1}\xi ,\ldots ,d_{x}^{k}\xi
\right) ,
\end{equation*}

where $d_{x}^{\ell }\xi $ stands for the $\ell ^{th}$ order differnetial at $%
x$ of section $\xi .$

Let $I$ be a $q$-tuple of non negative integers, say

\begin{equation*}
I=\left( i_{1},\ldots ,i_{q}\right)
\end{equation*}

with $0\leqslant i_{1},\ldots ,i_{q}\leqslant k.$ To every chain of order $%
\leqslant k,$ say $\theta \in C_{q}\left( \mathcal{G},W\right) ,$ and $%
\left( \xi _{1},\ldots ,\xi _{q}\right) \in \mathcal{G}^{q}$ we set

\begin{equation*}
\theta ^{I}\left( \xi _{1},\ldots ,\xi _{q}\right) \left( x\right) =\theta
\left( d_{x}^{i_{1}}\xi _{1},\ldots ,d_{x}^{i_{q}}\xi _{q}\right) .
\end{equation*}

Thus the chain $\theta $ is decomposed as it follows :

\begin{equation}
\theta =\sum_{I}\theta ^{I}  \tag{13}
\end{equation}

Then $\theta ^{I}$ is called component of type $I$ of $\theta .$

\bigskip

\textbf{Definition. }\textit{The symbol of a homogeneous }$q$-chain of order
$\leqslant k,$ say $\theta \in C_{q}\left( \mathcal{G},W\right) ,$ \textit{%
is its component of type} $\left( k,\ldots ,k\right) .$

\bigskip

The symbol of $\theta $ will be denoted by $\sigma _{\theta }.$

\bigskip \UNICODE{0x2026}

\textbf{Proposition. }\textit{The symbol }$\sigma _{\theta }$ o\textit{f
every homogeneous }$q$\textit{-cycle} $\theta $ \textit{is} $\delta _{q}$%
\textit{-closed.}

\bigskip

\textit{Proof.} Let us consider the decomposition of $\theta $ given by (13)

\begin{equation*}
\theta \left( \xi _{1},\ldots ,\right) \left( x\right) =\sum_{I}\theta
^{I}\left( \xi _{1},\ldots ,\xi _{q}\right)
\end{equation*}

with $I=\left( i_{1},\ldots ,i_{q}\right) .$ Then we see that

\begin{equation*}
\delta \theta \left( \xi _{1},\ldots ,\xi _{q+1}\right) =\sum_{I}\delta
\theta ^{I}\left( \xi _{1},\ldots ,\xi _{q+1}\right) .
\end{equation*}

If $q=0$ then $\theta \in J\left( W\right) ,$ since elements of $J\left(
W\right) $ are of order $\leqslant 0,$ we have $\sigma _{\theta }=0\quad
\forall \theta \in J\left( W\right) .$

Let us suppose that $q>0,$ then the formula

\begin{equation*}
\delta _{q}\theta \left( \xi _{1},\ldots ,\xi _{q+1}\right)
=\sum_{1\leqslant j\leqslant q}\left( -1\right) ^{j}\left\{ \left( \xi
_{j}\theta \right) \left( \ldots \hat{\xi}\ldots \xi _{q+1}\right) +\left(
e_{q}\left( \xi _{j}\right) \theta \xi _{q+1}\right) \left( \ldots \hat{\xi}%
,\ldots ,\hat{\xi}_{q+1}\right) \right\}
\end{equation*}

shows that the terms of $\delta \theta ^{I}$ are of order $\leqslant \max
\left( I\right) +1.$ Now we introduce the concept of degree of $\theta ^{I},$
say $\deg \theta ^{I}$ by setting

\begin{equation*}
\deg \theta ^{I}=\left| I\right| =i_{1}+\cdots +i_{q}
\end{equation*}

Thus if $\sigma _{\theta }\neq 0$ then $\deg \sigma _{\theta }=qk.$
Therefore we deduce that $\sigma _{\theta }\neq 0$ implies the inequality

\begin{equation*}
\deg \theta ^{I}<\deg \sigma _{\theta }.
\end{equation*}

On the other hand a simple calculation shows that $\delta \theta ^{I}\left(
\xi _{1},\ldots ,\xi _{q+1}\right) $ is homogeneous in the following sense :
let us set

\begin{equation}
\delta _{q}\theta ^{I}\left( \xi _{1},\ldots ,\xi q+1\right) =\sum_{J}\delta
_{q}\theta ^{I}\left( d^{j_{1}}\xi _{1},\ldots ,d^{j_{q+1}}\xi _{q+1}\right)
\tag{14}
\end{equation}

with $J=\left( j_{1},\ldots ,j_{q+1}\right) \in \mathbb{N}^{q+1}.$ Thus each
$J$ which occurs in $\left( 14\right) $ must satisfy the equality

\begin{equation*}
\left| J\right| =\left| I\right| +1
\end{equation*}

Thus if $\sigma _{\theta }\neq 0$ then $\theta ^{I}\neq \sigma _{\theta }$
implies that

\begin{equation}
\deg \left( \delta _{q}\theta ^{I}\right) <\deg \left( \delta _{q}\sigma
_{\theta }\right)  \tag{15}
\end{equation}

From $\left( 15\right) $ one concludes that $\delta _{q}\theta =0$ implies $%
\delta _{q}\sigma _{\theta }=0$ \ $\bullet $

\bigskip

\textbf{3.2. Symbols and transversally Poisson foliations.}

\qquad We are going to relate symbols of some kind of homogeneous 2-cycles
to foliations that have transverse Poisson structures.

To motivate we begin by recalling some elementary properties of Hochschild
complex of associative commutative algebras.

Let $\frak{A}$ be an associative commentative algebra over $\mathbb{F}.$ We
denote by $C^{2}\left( \frak{A,A}\right) $ the vector space of 2-cochains of
the Hochschild complex

\begin{equation*}
\frak{A}\overset{d}{\rightarrow }C^{1}\left( \frak{A,A}\right) \overset{d}{%
\rightarrow }\ldots C^{q}\left( \frak{A,A}\right) \overset{d}{\rightarrow }%
\ldots
\end{equation*}

The coboundary operator

\begin{equation*}
d:C^{2}\left( \frak{A,A}\right) \rightarrow C^{3}\left( \frak{A,A}\right)
\end{equation*}

is defined by

\begin{equation*}
d\theta \left( a,b,c\right) =a\theta \left( b,c\right) -\theta \left(
ab,c\right) +\theta \left( a,bc\right) -\theta \left( a,b\right) c.
\end{equation*}

We denote by $\Pi _{\theta }$ the skew symmetric component of $\theta ,$ viz

\begin{equation*}
2\Pi _{\theta }\left( a,b\right) =\theta \left( a,b\right) -\theta \left(
b,a\right)
\end{equation*}

The following elementary properties have remarkable consequences, [KM],
[BdM].

Suppose that $d\theta =0,$ then

$\left( P_{1}\right) $\quad $d\Pi _{\theta }=0\ ;$

$\left( P_{2}\right) \quad $for every fixed $a\in \frak{A}$ the linear map $%
\Pi _{\theta }\left( a,-\right) $ given by $b\rightarrow \Pi _{\theta
}\left( a,b\right) $ is a derivation of the algebra $\frak{A.}$

$\left( P_{3}\right) \quad $the symmetric component of $\theta ,$ say $%
S_{\theta }$ given by

\begin{equation*}
2S_{\theta }\left( a,b\right) =\theta \left( a,b\right) +\theta \left(
b,a\right)
\end{equation*}

is $d$-exact whenever $\frak{A}$ has a unit element.

Suppose that $\mathbb{F}$ is the field $\mathbb{R}$ of real numbers ; then
let $\frak{A}$ be the associative commutative algebra of smooth functions on
a smooth manifold $M.$ $\left( P_{1}\right) ,\left( P_{2}\right) $ and $%
\left( P_{3}\right) $ relate the isomorphism classes of star products on $M$
to isomorphism classes of formal Poisson structures (see [KM],[BdM],[VJ]).

An other reason why we are interessed in the scalar homology of
Koszul-Vinberg algebroid is the following. According to [BdM] it is unknown
whether does the correspondence

\begin{center}
\{class of star product on $\mathbb{N\}\leftrightarrow }$\{class of formal
Poisson structures on $M\}$
\end{center}

hold in the category of complex analytic manifolds.

\bigskip

\textbf{Remarks. }In general situation neither $\left( P_{1}\right) $ nor $%
\left( P_{2}\right) $ holds in Koszul-Vinberg algebras. For instance given
any Koszul-Vinberg algebra $\mathcal{A}$ its multiplication map

\begin{equation*}
\left( a,b\right) \rightarrow ab
\end{equation*}

is $\delta _{2}$-closed, (it is $\delta _{1}$-exact), however up to factor $%
\frac{1}{2}$ its skew symmetric component is nothing but the bracket of $%
\mathcal{A}_{L}$ ; that bracket is not $\delta _{2}$-closed.

Now in regard to $\left( P_{2}\right) $ the linear map $ad\left( a\right) $
will be a derivation of the Koszul-Vinberg algebra $\mathcal{A}$ iff $\xi
\in J\left( \mathcal{A}\right) .$

The remarks above motivate new definitions.

Given $\theta \in C_{2}\left( \mathcal{G},\mathcal{G}\right) $ we define the
$\theta $-associator to be the 3-multilinear map from $\mathcal{G}^{3}$ to $%
\mathcal{G}$ defined by

\begin{equation*}
\left( \xi _{1},\xi _{2},\xi _{3}\right) _{\theta }=\theta \left( \xi
_{1},\theta \left( \xi _{2},\xi _{3}\right) \right) -\theta \left( \theta
\left( \xi _{1},\xi _{2}\right) ,\xi _{3}\right)
\end{equation*}

\bigskip

\textbf{Definition.}\textit{\ (i) A chain} $\theta \in C_{2}\left( \mathcal{G%
},\mathcal{G}\right) $\textit{\ is called a Koszul-Vinberg chain of} $%
\mathcal{G}$ \textit{if} $\left( \xi _{1},\xi _{2},\xi _{3}\right) _{\theta
}=\left( \xi _{2},\xi _{1},\xi _{3}\right) _{\theta }$\quad $\forall \xi
_{1},\xi _{2},\xi _{3}\in \mathcal{G}\frak{.}$

\textit{(ii) A Koszul-Vinberg cycle of }$\mathcal{G}$\textit{\ is a
Koszul-Vinberg chain }$\theta $ \textit{satisfying both conditions} $\delta
_{2}\theta =0$ \textit{and} $\delta _{2}\Pi _{\theta }=0.$

\bigskip

Given a Koszul-Vinberg algebroid $\left( E,a_{n}\right) $ we consider its
complex

\begin{equation*}
J\left( W\right) \overset{\delta _{0}}{\rightarrow }C_{1}\left( \mathcal{G}%
,W\right) \overset{\delta _{1}}{\rightarrow }\ldots
\end{equation*}

Given a chain of order $\leqslant k,k>0,$ say $\theta \in C_{q}\left(
\mathcal{G},W\right) $ we decompose $\theta $ as before, that is

\begin{equation*}
\theta =\sum_{I}\theta ^{I}.
\end{equation*}

Henceforth let us suppose all the $\theta ^{I}$ to have the same degree, viz
$\left| I\right| =$constant$.$

\bigskip

\textbf{Theorem II. }\textit{Let} $\theta \in C_{0,2}\left( \frak{G}%
,W\right) \mathit{\ }$\textit{be a skew symmetric cycle of order} $\leqslant
k$ all of \textit{whose components} $\theta ^{I}$ \textit{have the same
degree equal to 2k. If }$k>0$ then $k=1.$

\bigskip

\textit{Proof.} We start by assuming that $k>1.$ Then let $f,g,h$ be
elements of $W=\Gamma \left( \mathbb{\tilde{F}}\right) .$ We consider the
following expression

\begin{equation*}
L\left( f,g,h\right) =f\theta \left( g,h\right) +\theta \left( g,f\right)
h-\theta \left( g,fh\right) .
\end{equation*}

From the closness of $\theta ,$ viz, $\delta _{2}\theta =0,$ we deduce that

\begin{equation*}
L\left( f,g,h\right) =L\left( g,f,h\right)
\end{equation*}

Since $\deg \theta =2k$ a direct calculation leads to

\begin{equation*}
L\left( f,g,h\right) =\sum_{\substack{ r+s=k  \\ rs>0}}\theta \left(
d^{k}g,d^{r}f\cdot d^{s}h\right)
\end{equation*}

where $d^{rf}\cdot d^{s}h\left( x\right) $ is the product of polynomial
functions. At the present step we know that $L\left( f,g,h\right) $ is
symmetric w.r.t. its three arguments $f,g$ and $h.$ This last property holds
iff $L\left( f,g,h\right) \equiv 0\quad \forall f,g,h\in W.$

Thus for every fixed $g\in W$ the linear map from $W$ to $W$ which is
defined by

\begin{equation*}
f\rightarrow \theta \left( g,f\right)
\end{equation*}

is a derivation of the associative algebra $W.$ So $\theta $ is a
bidifferential operator of order one, viz $k=1.$ This conclusion contredicts
our starting assumption $k>1.$ Theorem II is proved \ $\bullet $

\strut

The following statement is a direct consequence of Theorem II.

\bigskip

\textbf{Theorem III. }\textit{Let }$\left( E,a_{n}\right) $\textit{\ be a
Koszul-Vinberg algebroid. If }$\theta \in C_{0,2}\left( \mathcal{G},W\right)
$ \textit{is a Koszul-Vinberg cycle of order} $\leqslant k,$ \textit{then
the skew symmetric component of its symbol }$\sigma _{g}$ \textit{is Poisson
tensor on }$M.$

\bigskip

\textit{Proof.} We denote by $\Pi _{\theta }$ the skew symmetric component
of $\sigma _{\theta },$

\begin{equation*}
2\Pi _{\theta }\left( f,g\right) =\sigma _{\theta }\left( f,g\right) -\sigma
_{\theta }\left( g,f\right)
\end{equation*}

If $\Pi _{\theta }=0,$ then Theorem II is true.

Let us assume that $\Pi _{\theta }\neq 0,$ then all components of $\Pi
_{\theta },$ say $\Pi _{\theta }^{I},$ have degree $2k.$ By the virtue of
Theorem II we have $k=1.$ This last conclusion implies that

\begin{equation*}
\Pi _{\theta }\left( g,fh\right) =f\Pi _{\theta }\left( g,h\right) +\Pi
\left( g,f\right) h,\quad \forall f,g,h\in W.
\end{equation*}

Furthermore it is easy to see that $2\Pi _{\theta }$ is but the symbol of

\begin{equation*}
\Lambda _{\theta }\left( f,g\right) =\theta \left( f,g\right) -\theta \left(
g,f\right) .
\end{equation*}

Since $\theta $ is a Koszul-Vinberg chain one has

\begin{equation*}
\oint \Lambda _{\theta }\left( f,\Lambda _{\theta }\left( g,h\right) \right)
=0.
\end{equation*}

We observe that the 3-chain

\begin{equation*}
f,g,h\rightarrow \Lambda _{\theta }\left( f,\Lambda _{\theta }\left(
g,h\right) \right)
\end{equation*}

is of order $\leqslant 2.$

Regarding the degrees of components of $\Lambda _{\theta }\left( f,\Lambda
_{\theta }\left( g,h\right) \right) $ we conclude that the condition

\begin{equation*}
\oint \Lambda _{\theta }\left( f,\Lambda _{\theta }\left( g,h\right) \right)
=0
\end{equation*}

must implies the same condition on the symbol of $\Lambda _{\theta },$ viz

\begin{equation*}
\oint \Pi _{\theta }\left( f,\Pi _{\theta }\left( g,h\right) \right) =0
\end{equation*}

Now let us set

\begin{equation*}
\left\{ f,g\right\} =\Pi _{\theta }\left( f,g\right) ,\forall f,g\in W.
\end{equation*}

Then $\left( M,\left\{ ,\right\} \right) $ is Poisson manifold. This ends
the proof of Theorem III \ $\bullet $

\bigskip

Theorem III\textit{\ }is a useful tool to relate general 2-cycles $\theta
\in C_{2}\left( \mathcal{G},W\right) $\ to transverse Poisson structures of
Koszul-Vinberg algebroids.

\bigskip

We start with a Koszul-Vinberg algebroid $\left( E,a_{n}\right) $ and the
associated Koszul-Vinberg algebre $\mathcal{G}\frak{.}$

Let $\theta \in C_{2}\left( \mathcal{G},W\right) $ be a 2-cycle, viz $\delta
_{2}\theta =0.$ We shall decompose $\theta $ as follows :

\begin{equation*}
\theta =\theta _{2,0}+\theta _{1,1}+\theta _{0,2}\in C_{2,0}\left( \mathcal{G%
},W\right) \oplus C_{1,1}\left( \mathcal{G},W\right) \oplus C_{0,2}\left(
\mathcal{G},W\right) .
\end{equation*}

Given

\begin{equation*}
\xi =\left( s,f\right) ,\xi ^{\prime }=\left( s^{\prime },f^{\prime }\right)
\in \mathcal{G}\simeq \frak{\Gamma }\left( E\right) \oplus \Gamma \left(
\mathbb{\tilde{F}}\right) ,
\end{equation*}

we have

\begin{equation*}
\theta \left( \xi ,\xi ^{\prime }\right) =\theta _{2,0}\left( s,s^{\prime
}\right) +\theta _{1,1}\left( s,f^{\prime }\right) +\theta _{1,1}\left(
f,s^{\prime }\right) +\theta _{0,2}\left( f,f^{\prime }\right) .
\end{equation*}

By the virtue of Theorem I we know that $\delta _{2}\theta =0$ implies the
exactness of $\theta _{2,0}.$

To simplify the notations $s\cdot f$ will stand for $a_{n}\left( s\right)
f\quad \forall \left( s,f\right) \in \Gamma \left( E\right) \times \Gamma
\left( \mathbb{\tilde{F}}\right) .$

Thus given $s\in \mathcal{A=\Gamma }\left( E\right) $ and $f,g\in \mathcal{%
A=\Gamma }\left( \mathbb{\tilde{F}}\right) ,$\quad $\delta _{2}\theta =0$
implies that

\begin{eqnarray}
&&s\cdot \theta _{0,2}\left( f,g\right) -\theta _{0,2}\left( s\cdot
f,g\right) -\theta _{0,2}\left( f,s\cdot g\right) +\theta _{1,1}\left(
f,s\right) g  \TCItag{16} \\
&=&f\theta _{1,1}\left( s,g\right) -\theta _{1,1}\left( s,fg\right) -\theta
_{1,1}\left( s,f\right) g  \notag
\end{eqnarray}

and

\begin{equation}
f\theta _{1,1}\left( g,s\right) -\theta _{1,1}\left( fg,s\right) =g\theta
_{1,1}\left( f,s\right) -\theta _{1,1}\left( gf,s\right) .  \tag{17}
\end{equation}

From $\left( 17\right) $ we deduce that

\begin{equation*}
\theta _{1,1}\left( f,s\right) =f\theta _{1,1}\left( 1,s\right) \quad
\forall \left( s,f\right) \in \mathcal{A}\times W.
\end{equation*}

Then from (15) and (16) we must conclude that for any $\in \mathcal{A}$ and $%
f,g\in W$ we have

\begin{eqnarray}
s\theta _{0,2}\left( f,g\right) &=&\theta _{0,2}\left( s\cdot f,g\right)
+\theta _{0,2}\left( f,s\cdot g\right) -\theta _{1,1}\left( f,s\right) g
\TCItag{18} \\
&&+f\theta _{1,1}\left( s,g\right) -\theta _{1,1}\left( s,fg\right) +\theta
_{1,1}\left( s,f\right) g.  \notag
\end{eqnarray}

Let $\Lambda _{\theta }$ be the skew symmetric 2-chain in $C_{0,2}\left(
\frak{G},W\right) $ given by

\begin{equation*}
\Lambda _{\theta }\left( f,g\right) =\theta _{0,2}\left( f,g\right) -\theta
_{0,2}\left( g,f\right) .
\end{equation*}

From (18) we deduce that $\forall s\in \mathcal{A},$ $\forall f,g\in W$ one
has

\begin{eqnarray*}
s\Lambda _{\theta }\left( f,g\right) &=&\Lambda _{\theta }\left( s\cdot
f,g\right) +\Lambda _{\theta }\left( f,s\cdot g\right) \\
&&-\theta _{1,1}\left( f,s\right) g+\theta _{1,1}\left( g,g\right) f.
\end{eqnarray*}

We already know that $\theta _{1,1}\left( f,s\right) =f\theta _{1,1}\left(
1,s\right) .$ Thus (17) gives rise to the identity

\begin{equation}
s\Lambda _{\theta }\left( f,g\right) =\Lambda _{\theta }\left( s\cdot
f,g\right) +\Lambda _{\theta }\left( f,s\cdot g\right) .  \tag{19}
\end{equation}

We are now in position to prove the following statement

\bigskip

\textbf{Theorem IV.}\textit{\ Let }$\left( E,a_{n}\right) $ \textit{be a
regular Koszul-Vinberg algebroid. Then every Koszul-Vinberg 2-cycle of order
}$\leqslant k\mathit{\ }$\textit{with\ }$k>0$\textit{\ provides the
associative commutative algebra \ of first integrals of }$a_{n}\left(
E\right) $\textit{\ with a structure of Poisson algebra.}

\bigskip

\textit{Proof. }Let $\theta \in C_{2}\left( \mathcal{G},W\right) $ be a
Koszul-Vinberg cycle. We are concerned with its component $\theta _{0,2}.$
Let $\sigma _{\theta }$ be the symbol of $\theta _{0,2}$ and let $\Pi
_{\theta }$ be the skew symmetric part of $\sigma _{\theta }.$ If $\Pi
_{\theta }=0$ the conclusion of Theorem IV holds if we set

\begin{equation*}
\left\{ f,g\right\} =0\quad \forall f,g\in I\left( E\right) .
\end{equation*}

Let us suppose that $\Pi _{\theta }\neq 0.$ By the virtue of Theorem I we
have $k=1.$

Thus given $f,g\in W$ let us set

\begin{equation*}
\left\{ f,g\right\} =\Pi _{\theta }\left( f,g\right) .
\end{equation*}

According to Theorem III, $\Pi _{\theta }$ is a Poisson tensor on $M.$ Now
if $f$ and $g$ are elements of $I\left( E\right) ,$ then (19) shows that $%
\left\{ f,g\right\} $ is an element of $I\left( E\right) $ as well. That
ends the proof of Theorem IV \ $\bullet $

\strut

A particularly interesting situation is when a manifold is a locally flat
manifold, say $\left( M,D\right) .$ We consider the Koszul-Vinberg algebroid
$(TM,$ identity) ; recall that $TM$ stands for $T^{1,0}M$ when $\left(
M,D\right) $ is complex analytic. From our previous theorems above we deduce.

\bigskip

\textbf{Theorem V.}\textit{Given a locally flat manifold }$\left( M,D\right)
$\textit{\ every Koszul-Vinberg 2-cycle of order }$\leqslant k$\textit{\
with }$k>0$\textit{\ gives rise to a symplectic foliation in }$M$\textit{\
whose associative commutative algebra of first integrals is a trivial
Poisson algebra}.

\bigskip

\textit{Proof. }Theorem V is a particular case of Theorem IV. Indeed let $%
\theta \in C_{2}\left( \mathcal{G},W\right) $ be a Koszul-Vinberg 2-cycle, $%
\mathcal{G}$ being associated to the Koszul-Vinberg algebroid $(TM,$
identity). We assume $\theta $ to be of order $\leqslant k$ with $k>0.$ We
consider the skew symmetric part of $\sigma _{\theta }$ where $\sigma
_{\theta }$ is the symbol of $\theta _{0,2}$ ; thus we set

\begin{equation*}
2\Pi _{\theta }\left( f,g\right) =\sigma _{\theta }\left( f,g\right) -\sigma
_{\theta }\left( g,f\right) .
\end{equation*}

By Theorem IV we know that $\Pi _{\theta }$ is a Poisson tensor on $M.$ On
the other hand (18) shows that $\Pi _{\theta }$ is D-parallel, viz $\left(
M,D,\Pi _{\theta }\right) $ is affine Poisson manifold. That implies that
the differential system which is generated by hamiltonian vector fields is a
regular distribution on $M.$ So we obtain a foliation $\mathcal{F}_{\theta }$
whose tangent bundle is spanned by hamiltonian vector fields.

Of course each leaf of\/$\mathcal{F}_{\theta }$ is a symplectic manifold.
From (18) we deduce that the transverse Poisson tensor which is induced by $%
\Pi _{\theta }$ is zero. This ends the proof $\ \bullet $

\bigskip

\textit{Remark. }Theorem V can also be deduced from the local decomposition
theorem of A. Weinstein. [WA].

\bigskip

\textbf{4 - The cotangent Koszul-Vinberg algebroid.}

\bigskip

\qquad Let $M$ be a smooth manifold (viz $\mathbb{F=R)}.$

Recall that according to $\left[ \text{KM}\right] $ there is a one to one
correspondence between the set of isomorphism classes of star products on $M$
and the set of isomorphism classes of formal Poisson structures on $M.$ In
[BdM] Boutet de Monval point out that the arguments which are used to prove
the above equivalence don't have their analogues in the category of complex
analytic manifonds.The main aim of this section is to show that the concept
of the scalar homology of Koszul-Vinberg algebroid allows to treat
simultaneously the cases of smooth and complex analytic Poisson manifold
structures.

\bigskip

\textbf{4.1. The vertical Koszul-Vinberg algebroid.}

\qquad Let $M$ be a manifold and let $T^{\ast }M$ be its cotangent vector
bundle. When $M$ is complex analytic $T^{\ast }M$ stands for the vector
bundle whose section are holomorphic differential forms (viz $T^{\ast 1,0}M)$

The fibers of the canonical projection

\begin{equation*}
T^{\ast }M\overset{p}{\rightarrow }M
\end{equation*}

carry a canonical locally flat structure. They define in $T^{\ast }M$ a
foliation whose tangent bundle $\mathcal{V}\left( M\right) $ is a sub-bundle
of $T\left( T^{\ast }M\right) .$

The distribution $\mathcal{V}\left( M\right) \subset T\left( T^{\ast
}M\right) $ defines a lagrangian foliation in $\left( T^{\ast }M,\omega
_{\lambda }\right) $ where $\omega _{\lambda }$ is the Liouville symplectic
form on $T^{\ast }M$ ; so the structure of Koszul-Vinberg algebra of $%
\mathcal{A=\Gamma }\left( \mathcal{V}\left( M\right) \right) $ is that given
Example $\left( e_{2}\right) $ of sub-section 2.1.

We denote by Diff$_{0}\left( M\right) $ the group of diffeomorphisms of $M$
that are isotopic to identity map from $M$ to $M.$ Diffeomorphism of $M$
means smooth diffeomorphism (resp. complex analytic diffeomorphism) when $%
\mathbb{F=R}$ (resp. when $\mathbb{F=C)}$. Every $\varphi \in $ Diff$%
_{0}\left( M\right) $ gives rise to a $\mathcal{V}\left( M\right) $%
-preserving diffeomorphism of $T^{\ast }M,$ say $\varphi ^{\ast }.$

Actually $\varphi ^{\ast }$ is an automorphism of the Koszul-Vinberg algebra
$\mathcal{G}\frak{=}\mathcal{A\oplus \Gamma }\left( \mathbb{\tilde{F}}%
\right) $ where $\mathcal{A=\Gamma }\left( v\left( M\right) \right) $ and $%
\Gamma \left( \mathbb{\tilde{F}}\right) =\Gamma \left( T^{\ast }M\times
\mathbb{F}\right) .$

Henceforth we are interessed in the chain complexes

\begin{equation*}
J\left( W\right) \overset{\delta _{0}}{\rightarrow }C_{1}\left( \mathcal{G}%
,W\right) \overset{\delta _{1}}{\rightarrow }\cdots
\end{equation*}

and

\begin{equation*}
J\left( \mathcal{G}\right) \overset{\delta _{0}}{\rightarrow }C_{1}\left(
\mathcal{G},\mathcal{G}\right) \overset{\delta _{1}}{\rightarrow }\cdots
\end{equation*}

where $W=\Gamma \left( \mathbb{\tilde{F}}\right) .$

We are going to point out that $C\left( \mathcal{G},W\right) $ allows to
construct all Poisson structures on $M.$ Indeed let $\theta \in C_{2}\left(
\mathcal{G},W\right) $ be a 2-cycle of order $\leqslant k$ with $k>0.$

To every such a cycle, say $\theta $, we assign the skew symmetric part of $%
\theta _{0,2}.$

Thus we set

\begin{equation*}
\Lambda _{\theta }\left( f,g\right) =\theta _{0,2}\left( f,g\right) -\theta
_{0,2}\left( g,f\right) \quad \forall f,g\in W.
\end{equation*}

We know that the symbol $\Pi _{\theta }$ of $\Lambda _{\theta }$ is but the
skew part of the symbol of $\theta _{0,2}.$

Suppose that $\Pi _{\theta }$ is a Poisson tensor on $T^{\ast }M,$ then it
is the unique Poisson tensor in its scalar homology class

\begin{equation*}
\left[ \Pi _{\theta }\right] \in H_{0,2}\left( v\left( M\right) ,\mathbb{F}%
\right) \text{.}
\end{equation*}

The following definition is useful to our purpose.

\bigskip

\textbf{Definition. }\textit{A homology class }$c\in H_{2}\left( v\left(
M\right) ,\mathbb{F}\right) $\textit{\ is called a Poisson class if it
contains a cycle of order }$\leqslant 1,$\textit{\ say }$\theta $\textit{\
such that }$\Pi _{\theta }$\textit{\ is a Poisson tensor on }$T^{\ast }M.$

\bigskip

Let $\mathbb{F=R}$ and $\frak{A}$ be the associative commutative algebra $%
C^{\infty }\left( M,\mathbb{R}\right) $of smooth functions on a smooth
manifold $M.$

Then considering the Hochsild complex of $\frak{A}$

\begin{equation*}
\frak{A}\overset{d}{\rightarrow }C^{1}\left( \frak{A},\frak{A}\right)
\overset{d}{\rightarrow }\cdots
\end{equation*}

every class $c\in H^{2}\left( \frak{A},\frak{A}\right) $ is a Poisson class.

According to our previous definitions every homology class $c$ which
contains a Koszul-Vinberg cycle is a Poisson class.

The set of Poisson classes in $H_{2}\left( v\left( M\right) ,\mathbb{F}%
\right) $ is denoted by $\mathcal{P}H_{2}\left( \mathbb{F}\right) $ and the
set of Poisson structures on $M$ is denoted by $\mathcal{P}\left( M\right) .$

Given a Poisson class $c\in H\left( v\left( M\right) ,\mathbb{F}\right) $
containing a 2-cycle $\theta $ such that $\Pi _{\theta }$ is a Poisson
tensor on $T^{\ast }M,$ let $f$ and $g$ be two first integrals of $v\left(
M\right) ,$ then both $f$ and $g$ have the form

\begin{eqnarray*}
f &=&f_{0}\circ p, \\
g &=&g_{0}\circ p
\end{eqnarray*}

where $f_{0}$ and g$_{0}$ are elements of $\Gamma \left( M\times \mathbb{F}%
\right) .$

By the virtue of (19) we have

\begin{equation}
s\cdot \Pi _{\theta }\left( f,g\right) =0\quad \forall s\in \mathcal{%
A=\Gamma }\left( v\left( M\right) \right) .  \tag{20}
\end{equation}

Thus the manifold $M$ carries a Poisson structure $\left( M,\left\{
{}\right\} _{\theta }\right) $ which is defined by

\begin{equation}
\left\{ f_{0},g_{0}\right\} =\Pi _{\theta }\left( f_{0}\circ p,g_{0}\circ
p\right)  \tag{21}
\end{equation}

where $p$ is the canonical projection from $T^{\ast }M$ to $M.$ Since a
class $c\in H_{2}\left( v\left( M\right) ,\mathbb{F}\right) $ contains at
most one Poisson tensor the map $\mathcal{PH}_{2}\left( \mathbb{F}\right)
\overset{\Pi }{\rightarrow }\mathcal{P}\left( M\right) $ is well defined.

\bigskip

\textbf{Theorem VI. }\textit{Given a manifold }$M$\textit{\ the map }$\Pi $%
\textit{\ which assigns to every class }$\left[ \theta \right] =c\in
\mathcal{P}H_{2}\left( \mathbb{F}\right) $\textit{\ the Poisson tensor }$\Pi
_{\theta }\in \mathcal{P}\left( M\right) $\textit{\ is an isomorphism on }$%
\mathcal{P}\left( M\right) .$

\bigskip

\textit{Proof.} We have but to prove that $\Pi :\mathcal{PH}_{2}\left(
\mathbb{F}\right) \rightarrow \mathcal{P}\left( M\right) $ is onto. Let $P$
be a Poisson tensor on $M.$We denote by

\begin{equation*}
0:M\rightarrow T^{\ast }M
\end{equation*}

the zero section of $p:T^{\ast }M\rightarrow M.$ The induced linear map

\begin{equation*}
\sigma ^{\ast }:\text{ Funct }\left( T^{\ast }M,\mathbb{F}\right)
\rightarrow \text{ Funct}\left( M,\mathbb{F}\right)
\end{equation*}

is surjective. Funct$\left( M,\mathbb{F}\right) $ is the associative
commutative algebra $\mathbb{F}$-valued functions on $M.$ Of course Funct$%
\left( M,\mathbb{F}\right) $ is the algebra of smooth functions (resp
holomorphic functions) when $\mathbb{F=R}\left( \text{resp.}\mathbb{F=C}%
\right) .$Let us define the bracket $\tilde{P}$\ on Funct$\left( T^{\ast }M,%
\mathbb{F}\right) $ by setting

\begin{equation}
\tilde{P}\left( f,g\right) =p^{\ast }P\left( \sigma ^{\ast }\left( f\right)
,v^{\ast }\left( g\right) \right) \quad \forall f,g\in Funct\left( T^{\ast
}M,\mathbb{F}\right) .  \tag{22}
\end{equation}

Since both

\begin{eqnarray*}
p^{\ast } &:&\text{Funct }\left( M,\mathbb{F}\right) \rightarrow \text{Funct
}\left( T^{\ast }M\right) \text{ and} \\
\sigma ^{\ast } &:&\text{Funct }\left( T^{\ast }M,\mathbb{F}\right)
\rightarrow \text{Funct }\left( M,\mathbb{F}\right)
\end{eqnarray*}

are homomorphisms of associative algebras with $v^{\ast }\circ p^{\ast
}\left( f\right) =f$ we have

\begin{equation*}
\oint \tilde{P}\left( f,\tilde{P}\left( g,h\right) \right) =0
\end{equation*}

and

\begin{equation}
\tilde{P}\left( g,fh\right) =f\tilde{P}\left( gh\right) +\tilde{P}\left(
g,f\right) h  \tag{23}
\end{equation}

Thus $\tilde{P}$ is a Poisson tensor on $T^{\ast }M.$

It remains to see that $\tilde{P}\in C_{0,2}\left( \mathcal{G}\text{,W}%
\right) $ is $\delta _{2}$-closed .We fix a system of local coordinate
functions on $M,$ say $q=\left( q_{1}\cdots q_{m}\right) ,$ whose domain is
an open set $U\subset M$ . Then $p^{-1}\left( U\right) =T^{\ast }U$ carry
the local coordinate functions

\begin{equation*}
\left( q,p\right) =\left( q_{1}\cdots q_{m},p_{1}\cdots p_{m}\right) .
\end{equation*}

So $\forall \xi \in T^{\ast }U$ we have

\begin{equation*}
\xi =\sum_{j=1}^{m}p_{j}dq_{j}
\end{equation*}

Using those coordinate functions we have

\begin{equation}
\tilde{P}\left( \xi \right) =\sum_{i,j}P_{ij}\left( q_{1}\cdots q_{m}\right)
\frac{\partial }{\partial p_{i}}\wedge \frac{\partial }{\partial qj}
\tag{24}
\end{equation}

Every $S\in \mathcal{A=\Gamma }\left( v\left( M\right) \right) $ has the form

\begin{equation*}
S\left( \xi \right) =\sum_{i=1}^{m}S_{i}\left( q,p\right) \frac{\partial }{%
\partial p_{i}}
\end{equation*}

Therefore we see that

\begin{equation}
S\cdot \tilde{P}\left( f,g\right) -\tilde{P}\left( s\cdot f,g\right) -\tilde{%
P}\left( f,s\cdot g\right) =0  \tag{25}
\end{equation}

for arbitrary elements $f,g\in W=$ Funct $\left( T^{\ast }M,\mathbb{F}%
\right) .$ From (23) and (25) we conclude that

\begin{equation*}
\delta _{2}\tilde{P}=0.
\end{equation*}

To end the proof we define $\theta \in C_{2}\left( \mathcal{G},W\right) $ by
putting

\begin{equation*}
\theta \left( \left( s,f\right) ,\left( s^{\prime },f^{\prime }\right)
\right) =\frac{1}{2}\tilde{P}\left( f,f^{\prime }\right) .
\end{equation*}

Then $\left( 23\right) $ and $\left( 25\right) $ imply that

\begin{equation*}
\delta _{2}\theta =0
\end{equation*}

Actually we have

\begin{equation*}
\Pi _{\theta }=\tilde{P}
\end{equation*}

Therefore $\left[ \theta \right] \in \mathcal{PH}_{2}\left( \mathbb{F}%
\right) $ and we conclude that $\Pi $ is surjective $\ \bullet $

\strut

For $\mathbb{F=R}$ the homology space $H_{1}\left( \mathcal{V}\left(
M\right) \right) $ is of infinite dimension. That is a direct corollary of
the following statement.

\bigskip

\textbf{Theorem VII. }\textit{Let }$M$\textit{\ be a smooth manifold, then
there is a canonical linear injection of }$\Gamma \left( TM\right) $\textit{%
\ in }$H_{1}\left( \mathcal{V}\left( M\right) \right) .$

\bigskip

\textit{Proof. }Let X$\in \Gamma \left( TM\right) .$ Every $x\in M$ have a
neighbourhood $U\left( x\right) $ on which $X$ defines a local one parameter
group $\varphi _{t}\in $ Diff$_{0}\left( U\left( x\right) \right) ,$ $\left|
t\right| <\alpha .$ Then $\varphi _{t}^{\ast }$ is a local $v\left( M\right)
$-preserving diffeomorphism of $T^{\ast }\left( M\right) .$ Thus we define a
local section of $TT^{\ast }M,$ say $\tilde{X}$ by setting

\begin{equation*}
\tilde{X}_{u}\left( \xi \right) =\frac{d}{dt}\varphi _{t}\left( \xi \right)
_{\left| t=0\right. },\forall \xi \in T^{\ast }U\left( x\right) .
\end{equation*}

Actually $\tilde{X}$ is defined globally because on $T^{\ast }U\left(
x\right) \cap T^{\ast }U\left( y\right) $ we have $\tilde{X}_{U\left(
x\right) }=\tilde{X}_{U\left( y\right) }.$

On the other hand $\varphi _{t}^{\ast }$ is a local automorphism of the
Koszul-Vinberg algebra $\mathcal{G}\frak{.}$Thus the Lie derivation $w\cdot
r\cdot t$ $\tilde{X}$ is a derivation of the algebra $\mathcal{G}\frak{.}$
(That statement may be deduced from the fact that $\varphi _{t}^{\ast }$
preserves both $v\left( M\right) $ and the symplectic form $\omega _{\lambda
}.)$ We denote the vector space spanned by the $\tilde{X}^{^{\prime }s}$ by $%
\mathcal{\tilde{X}}\left( M\right) .$ We just saw that every Lie derivation $%
L_{\tilde{X}}$ is 1-cycle in $C_{1}\left( \mathcal{G},\mathcal{G}\frak{.}%
\right) .$

Since

\begin{equation*}
J\left( \mathcal{G}\frak{.}\right) \cap \mathcal{\tilde{X}}\left( M\right) =0
\end{equation*}

we see that $\forall \tilde{X}\in \mathcal{\tilde{X}}\left( M\right) $ the
homology class

\begin{equation*}
\left[ L_{\tilde{X}}\right] \in H_{1}\left( v\left( M\right) \right)
\end{equation*}

is different from zero. So the map

\begin{equation*}
X\in \mathcal{X}\left( M\right) \rightarrow \left[ L_{\tilde{X}}\right] \in
H_{1}\left( v\left( M\right) \right)
\end{equation*}

is injective. That ends the proof of Theorem VII \ $\bullet $

\bigskip

\textbf{5 - Poisson tensors as contact }$KV$\textbf{-invariants}

\bigskip

\qquad Section 4 provides a universal way to get all Poisson structures on
smooth (resp. complex analytic) manifolds. The aim of this section is to
investigate relation-ships between Poisson tensors and contact $KV$%
-invariants.

\bigskip

\textbf{5.1. Koszul-Vinberg algebroid of contact structure.}

\qquad Let $\left( M^{2n+1},\alpha \right) $ be a smooth contact manifold.
Then the Reeb vector field $\mathcal{R}$ defines canonically a 1-rank
sub-bundle of $TM$, say $E_{\mathcal{R}}$. Sections of $E_{\mathcal{R}}$ are
vector fields $f\mathcal{R}$ with $f\in W=\Gamma \left( \mathbb{\tilde{R}}%
\right) .$

We define multiplication in $\Gamma \left( E_{\mathcal{R}}\right) $ by

\begin{equation*}
f\mathcal{R}\cdot g\mathcal{R}=f<dg,\mathcal{R}>\mathcal{R}.
\end{equation*}

It easy to verify that

\begin{equation*}
\left( f\mathcal{R},g\mathcal{R},h\mathcal{R}\right) =\left( g\mathcal{R},f%
\mathcal{R},h\mathcal{R}\right) ,\forall f,g,h\in \Gamma \left( \mathbb{%
\tilde{R}}\right) .
\end{equation*}

We define the anchor $a_{n}$ by

\begin{equation*}
a_{n}\left( f\mathcal{R}\right) =f\mathcal{R}.
\end{equation*}

Thus $\left( E_{\mathcal{R}},\text{identity}\right) $ is a regular
Koszul-Vinberg algebra.

Let $\varphi \in $Diff$_{\alpha }\left( M\right) $ be a $\alpha $-preserving
diffeomorphism (contact diffeomorphism). Then $\varphi $ is $\mathcal{R}$%
-preserving, viz we have

\begin{equation*}
\text{Diff}_{\alpha }\left( M\right) =\text{Diff}_{\alpha ,\mathcal{R}%
}\left( M\right) .
\end{equation*}

Thus by setting $\varphi ^{\ast }h=h\circ \varphi ^{-1}\quad \forall h\in W,$%
one has

\begin{eqnarray}
\varphi _{\ast }\left( f\mathcal{R}\cdot g\mathcal{R}\right) &=&\varphi
_{\ast }\left( f\mathcal{R}\right) \cdot \varphi _{\ast }\left( g\mathcal{R}%
\right) .  \TCItag{26} \\
\varphi ^{\ast }\left( \mathcal{R}h\right) &=&\varphi _{\ast }\mathcal{R}%
\left( \varphi ^{\ast }h\right) .  \notag
\end{eqnarray}

Henceforth we are concerned with the Koszul-Vinberg algebra

\begin{equation*}
\mathcal{G}=\Gamma \left( E_{\mathcal{R}}\right) \oplus \mathbb{\tilde{R}}
\end{equation*}

whose multiplication is given by (9). Identities (26) show that every $%
\varphi \in $Diff$_{\alpha }\left( M\right) $ gives rise to an automorphism
of the complexes

\begin{eqnarray}
&&J\left( W\right) \overset{\delta _{0}}{\rightarrow }C_{1}\left( \mathcal{G}%
,W\right) \overset{\delta _{1}}{\rightarrow }\cdots  \TCItag{27} \\
&&J\left( \mathcal{G}\right) \overset{\delta _{0}}{\rightarrow }C_{1}\left(
\mathcal{G},\mathcal{G}\right) \overset{\delta _{1}}{\rightarrow }\cdots .
\notag
\end{eqnarray}

According to our previous notations, the group Diff$_{\alpha }\left(
M\right) $ acts on both $H\left( E_{\mathcal{R}},\mathbb{R}\right) $ and $%
H\left( E_{\mathcal{R}}\right) ,$ where

\begin{eqnarray*}
H\left( E_{\mathcal{R}},\mathbb{R}\right) &=&\bigoplus\limits_{q}H_{q}\left(
E_{\mathcal{R}},\mathbb{R}\right) \\
H\left( E_{\mathcal{R}}\right) &=&\bigoplus\limits_{q}H_{q}\left( E_{%
\mathcal{R}}\right) .
\end{eqnarray*}

Thus the homology spaces $H\left( E_{\mathcal{R}},\mathbb{R}\right) $ and $%
H\left( E_{\mathcal{R}}\right) $ yield new invariants of the contact
structure $\left( M,\alpha \right) .$

We also remark that from the actions of Diff$_{\alpha }\left( M\right) $ on
both $C\left( \mathcal{G},W\right) $ and $C\left( \mathcal{G},\mathcal{G}%
\right) $ we can deduce the equivariant subcomplexes

\begin{eqnarray}
&&J\left( W\right) ^{G}\overset{\delta _{0}}{\rightarrow }C_{1}^{G}\left(
\mathcal{G},\mathcal{W}\right) \overset{\delta _{1}}{\rightarrow }\cdots
\TCItag{28} \\
&&J\left( \mathcal{G}\right) ^{G}\overset{\delta _{0}}{\rightarrow }%
C_{1}^{G}\left( \mathcal{G},\mathcal{G}\right) \overset{\delta _{1}}{%
\rightarrow }\cdots .  \notag
\end{eqnarray}

In (28) $G=$Diff$_{\alpha }\left( M\right) ,$ $C_{q}^{G}\left( \mathcal{G},%
\mathcal{W}\right) $ and $C_{q}^{G}\left( \mathcal{G},\mathcal{G}\right) $
are the vector space of $G$-equivariant $q$-chains. The last notion makes
sense because $G$ acts on both $\mathcal{G}$ and $W.$ So, given $\theta \in
C_{q}\left( \mathcal{G},W\right) $ then $\theta \in C_{q}^{G}\left( \mathcal{%
G},W\right) $ iff we have the identity

\begin{equation*}
\varphi ^{\ast }\left( \theta \left( \xi _{1}\cdots \xi _{q}\right) \right)
=\theta \left( \tilde{\varphi}\left( \xi _{1}\right) ,\ldots ,\tilde{\varphi}%
\left( \xi _{q}\right) \right) .
\end{equation*}

Let us recall that $\forall \xi =\left( s,f\right) \in \Gamma \left( E_{%
\mathcal{R}}\right) \oplus W,\quad \forall \varphi \in G,$ then

\begin{equation}
\tilde{\varphi}\left( \xi \right) =\left( \varphi _{\ast }S,f\circ \varphi
^{-1}\right)  \tag{29}
\end{equation}

The $q^{th}$ homology spaces of $C^{G}\left( \mathcal{G}\text{,W}\right) $
and $C^{G}\left( \mathcal{G}\text{,}\mathcal{G}\right) $ are denoted by $%
H_{q}^{e}\left( E_{\mathcal{R}},\mathbb{R}\right) $ and $H_{q}^{e}\left( E_{%
\mathcal{R}}\right) $ respectively.

Classically, the vector spaces of residual cycles are denoted by $%
Z_{q}^{r}\left( \mathcal{G},W\right) \subset C_{q}\left( \mathcal{G}%
,W\right) $ and $Z_{q}^{r}\left( \mathcal{G},\mathcal{G}\right) \subset
C_{q}\left( \mathcal{G},\mathcal{G}\right) .$ An element $\theta \in
Z_{q}^{r}\left( \mathcal{G},W\right) $ iff $\delta _{q}\theta $ is $G$%
-equivariant. The vector space of residual boundaries is defined to be

\begin{equation*}
B_{q}^{r}=\varphi _{q-1}\left( C_{q-1}\left( \mathcal{G},W\right)
+C_{q}^{G}\left( \mathcal{G},W\right) \subset Z_{q}^{r}\left( \mathcal{G}%
,W\right) \right) .
\end{equation*}

It is well known that by setting

\begin{equation*}
H_{q}^{r}\left( E_{\mathcal{R}},\mathbb{R}\right) =Z_{q}^{r}\left( \mathcal{G%
},W\right) /B_{q}^{r}\left( \mathcal{G},W\right)
\end{equation*}

We obtain the following exact sequences

\begin{equation}
\cdots \rightarrow H_{q}^{e}\left( E_{\mathcal{R}},\mathbb{R}\right)
\rightarrow H_{q}\left( E_{\mathcal{R}},\mathbb{R}\right) \rightarrow
H_{q}^{r}\left( E_{\mathcal{R}},\mathbb{R}\right) \rightarrow
H_{q+1}^{e}\left( E_{\mathcal{R}},\mathbb{R}\right) \rightarrow \cdots
\tag{30}
\end{equation}

and

\begin{equation}
\cdots \rightarrow H_{q}^{e}\left( E_{\mathcal{R}}\right) \rightarrow
H_{q}\left( E_{\mathcal{R}}\right) \rightarrow H_{q}^{r}\left( E_{\mathcal{R}%
}\right) \rightarrow H_{q+1}^{e}\left( E_{\mathcal{R}}\right) \rightarrow
\cdots  \tag{31}
\end{equation}

of course the group $G=$Diff$_{a}\left( M\right) $ preserves the
bigraduations

\begin{eqnarray*}
C_{q}\left( \mathcal{G},W\right) &=&\bigoplus_{r+s=q}C_{r,s.}\left( \mathcal{%
G},W\right) , \\
C_{q}\left( \mathcal{G},\mathcal{G}\right) &=&\bigoplus_{r+s=q}C_{r,s}\left(
\mathcal{G},\mathcal{G}\right) .
\end{eqnarray*}

We observe that not all of the vector spaces $C_{q}^{G}\left( \mathcal{G}%
,W\right) $ are trivial.

Indeed let us consider the volume form

\begin{equation*}
v=\alpha \wedge \left( d\alpha \right) ^{n.}
\end{equation*}

It defines an isomorphism $\beta $ from $\Gamma \left(
\bigwedge\limits^{k}TM\right) $ to $\Gamma \left(
\bigwedge\limits^{2n+1-k}T^{\ast }M\right) $ for $0\leqslant k\leqslant
2n+1. $ For $\xi \in \Gamma \left( \bigwedge\limits^{k}TM\right) $ we set

\begin{equation*}
\beta \left( \xi \right) =i\left( \xi \right) v
\end{equation*}

where $i\left( \xi \right) $ is the inner product by $\xi .$ The isomorphism
$\beta $ is $G$-equivariant, that is to say

\begin{equation}
\beta \left( \varphi _{\ast }\xi \right) =\varphi ^{\ast }\left( \beta
\left( \xi \right) \right) ,\forall \varphi \in G.  \tag{32}
\end{equation}

The action of $G$ on $\Gamma \left( \Lambda T^{\ast }M\right) $ is given by $%
\left( 27\right) .$

Regarding $\tau ^{\ast }\left( M\right) =\bigoplus\limits_{r}\Gamma \left(
\bigotimes\limits^{r}T^{\ast }M\right) $ as $\mathbb{R}$-algebra, we denote
by $\tau _{\alpha }^{\ast }\left( M\right) $ the real subalgebra of $\tau
^{\ast }\left( M\right) $ that is generated by $\alpha $ and $d\alpha .$ Let
us set

\begin{equation*}
\#=\left( \beta \right) ^{-1}
\end{equation*}

Then $\#$ is an isomorphism from $\tau ^{\ast }\left( M\right) $ onto $\tau
\left( M\right) $ where

\begin{equation*}
\tau \left( M\right) =\bigoplus_{r}\Gamma \left( \bigotimes^{r}TM\right) .
\end{equation*}

According to (3) we get $\tau _{\alpha }\left( M\right) =\#\left( T_{\alpha
}^{\ast }\left( M\right) \right) $

We have in particular that

\begin{equation*}
\#\left( \left( d_{\alpha }\right) ^{n}\right) =\mathcal{R}
\end{equation*}

where $\mathcal{R}$ is the Reeb vector field of $\left( M,\alpha \right) .$
Every element $\xi \in \tau _{\alpha }\left( M\right) $ is $G$-invariant,
that is a straight consequence of (32).

Now we may regard elements of $\tau _{\alpha }\left( M\right) $ as operators
on $W=C^{\infty }\left( M,\mathbb{R}\right) ,$ so homogeneous element of
degre $m\in \mathbb{N},\left( \deg \left( \alpha \right) =1\text{ and }\deg
\left( d\alpha \right) =2\right) ,$ will be regarded as elements of $%
C_{o,m}^{G}\left( \mathcal{G},W\right) .$ This proves that the complex

\begin{equation*}
J\left( W\right) ^{G}\overset{\delta _{0}}{\rightarrow }C_{1}^{G}\left(
\mathcal{G},W\right) \overset{\delta _{1}}{\rightarrow }\cdots
\end{equation*}

has non zero chains.

To motivate we will consider the sub-complex

\begin{equation}
J\left( W\right) \overset{\delta _{0}}{\rightarrow }C_{1}\left( \mathcal{G}_{%
\mathbb{R}},W\right) \overset{\delta _{1}}{\rightarrow }\cdots  \tag{33}
\end{equation}

where

\begin{equation*}
\mathcal{G}_{\mathbb{R}}=\mathbb{R\cdot }\mathcal{R+}W
\end{equation*}

$\mathcal{G}_{\mathbb{R}}$ is a subalgebra of the Koszul-Vinberg algebra

\begin{equation*}
\mathcal{G=C}^{\infty }\left( M,\mathbb{R}\right) \cdot \mathcal{R+}W.
\end{equation*}

and it makes sense to consider (33) ; we have the following results.

Considering the vector space $W$ as a Koszul-Vinberg module of $\mathcal{G}_{%
\mathbb{R}}$

$J\left( W\right) $ is but the vector space of first integrals of $\mathcal{R%
}$, say $I\left( \mathcal{R}\right) .$ Thus in the chain complex (33), that
is

\begin{equation*}
I\left( W\right) \overset{\delta _{0}}{\rightarrow }C_{1}\left( \mathcal{G}_{%
\mathbb{R}},W\right) \overset{\delta _{1}}{\rightarrow }\cdots
\end{equation*}

we have $\delta _{0}=0$ and $H_{0}\left( \mathcal{G}_{\mathbb{R}},W\right)
=I\left( \mathcal{R}\right) .$

Now the complex (33) is stable under the action of $G.$ Not all $G$%
-invariant chains of $Cq\left( \mathcal{G}_{\mathbb{R}},W\right) $ are zero.
Indeed using the isomorphism \# we see for every non negative integer $%
m\quad 0\leqslant m\leqslant n$ the $m$-multivector,

\begin{equation*}
\Pi _{m}=\#\left( \alpha \wedge \left( d\alpha \right) ^{m}\right)
\end{equation*}

is $G$-invariant. Therefore we see that the chains

\begin{equation*}
\alpha ^{\ell }\otimes \Pi _{m}\in C_{\ell ,2n+1-m}\left( \mathcal{G}_{%
\mathbb{R}},W\right) ,\ell \in \mathbb{N}
\end{equation*}

are $G$-equivariant $\alpha ^{\ell }$ stands for $\bigotimes\limits^{\ell
}\alpha .$

So (33) contains the (non trivial) subcomplex

\begin{equation}
I\left( \mathcal{R}\right) ^{G}\overset{\delta _{0}}{\rightarrow }%
C_{1}^{G}\left( \mathcal{G}_{\mathbb{R}},W\right) \overset{\delta _{1}}{%
\rightarrow }\cdots C_{q}^{G}\left( \mathcal{G}_{\mathbb{R}},W\right)
\overset{\delta _{1}}{\rightarrow }\cdots  \tag{34}
\end{equation}

We consider $\mathcal{R}$ as element of $C_{0,1}^{G}\left( \mathcal{G}_{%
\mathbb{R}},W\right) ,$ thus one has

\begin{equation*}
\delta _{1}\mathcal{R}=0.
\end{equation*}
\quad

(Moreover $\mathcal{R}$ is of order $\leqslant 1$ and $\sigma _{\mathcal{R}}=%
\mathcal{R)}$

Since $\delta _{0}=0$ we the homology space $\tilde{H}_{1}\left( \mathcal{G}%
_{\mathbb{R}},W\right) $ of (34) is non zero. An other example of non zero
homology class of (34) is

\begin{equation*}
\alpha \otimes \mathcal{R}\in C_{1,1}^{G}\left( \mathcal{G}_{\mathbb{R}%
},W\right) .
\end{equation*}

We have $\delta _{1}\left( \alpha \otimes \mathcal{R}\right) =0$ and $\left[
\alpha \otimes \mathcal{R}\right] \in \tilde{H}_{2}\left( \mathcal{G}_{%
\mathbb{R}},W\right) \mathbb{-}\left\{ 0\right\} .$

Actually one can use the isomorphism \# to produce other canonical $G$%
-equivarianrt chaim in $C_{q}\left( \mathcal{G}_{\mathbb{R}},W\right) .$
Indeed let $\left( \ell ,m\right) \in \mathbb{N\times N}$ with $0\leqslant
m\leqslant n.$ We set

\begin{equation*}
\Pi _{2\left( n-m\right) }=\#\left( \alpha \wedge \left( d\alpha \right)
^{m}\right)
\end{equation*}

Thus $\Pi _{2\left( n-m\right) }$ is a G-invariant section of $%
\bigwedge\limits^{2\left( n-m\right) }TM.$We regard

\begin{equation*}
\theta =\alpha ^{\ell }\otimes \Pi _{2\left( n-m\right) }
\end{equation*}

as an element of $C_{\ell ,2\left( n-m\right) }^{G}\left( \mathcal{G}_{%
\mathcal{R}},W\right) $.

Actually every homology class of the complex

\begin{equation*}
I\left( \mathcal{R}\right) ^{G}\overset{\delta _{0}}{\rightarrow }%
C_{1}^{G}\left( \mathcal{G}_{\mathbb{R}},W\right) \overset{\delta _{1}}{%
\rightarrow }\cdots \rightarrow C_{q}^{G}\left( \mathcal{G}_{\mathbb{R}%
},W\right) \overset{\delta _{q}}{\rightarrow }\cdots
\end{equation*}

is a contact invariant of $\left( M,\alpha \right) .$

Remark : one easily checks that

\begin{equation*}
L_{\mathcal{R}}\Pi _{2\left( n-m\right) }=0
\end{equation*}

In particular $\Pi _{2}=\#\left( \alpha \wedge \left( d\alpha \right)
^{n-1}\right) $ is such that

\begin{equation*}
L_{\mathcal{R}}\Pi _{L}=0
\end{equation*}

we consider the foliation of $M$ which is defined by $\mathcal{R}.$

Thus given two functions $f_{1}$ and $f_{2}$ we have

\begin{equation}
\delta _{2}\Pi _{2}\left( \mathcal{R},f_{1},f_{2}\right) =0  \tag{35}
\end{equation}

Actually it is an exercice to show that $\delta _{2}\Pi _{2}=0$ and that $%
\Pi _{2}$ is but the Poisson tensor that is associated to the transverse
symplectic form $d\alpha .$

A Poisson tensor $\Pi _{2}\in C_{0,2}\left( \mathcal{G}_{\mathcal{R}%
},W\right) $ cannot be homologuous to zero. Since $C_{0,2}\left( \mathcal{G}%
_{\mathbb{R}},W\right) =C_{0,2}\left( \mathcal{G},W\right) $ we deduce from
(35) that $\Pi _{2}$ is a G-invariant cycle in $C_{0,2}\left( \mathcal{G}%
,W\right) $. Thus we have

\begin{equation*}
\left[ \Pi _{2}\right] \in H_{2}^{L}C_{0,2}C_{0,2}\left( \mathcal{G}%
,W\right) -\left\{ 0\right\}
\end{equation*}

More generally we

\begin{eqnarray*}
\Pi _{2\left( n-m\right) } &\in &C_{0,2\left( n-m\right) }\left( \mathcal{G}%
,W\right) ^{G} \\
\Pi _{2\left( n-m\right) +1} &=&\#\left( d\alpha \right) ^{m}\in
C_{0,2\left( n-m\right) +1}^{G}\left( \mathcal{G},\mathcal{W}\right) .
\end{eqnarray*}

\begin{center}
\bigskip \bigskip \strut \bigskip \strut

\textbf{Bibliographie}

\bigskip

\bigskip
\end{center}

$\left[ \text{\textbf{BdM}}\right] $ \ \textbf{BOUTET \ de \ MONVEL \ L. }%
Complex star algebras. Math. Phys. An. Geom. 2 (1999) n$%
{{}^\circ}%
$2 113-139.

\bigskip

$\left[ \text{\textbf{DLM}}\right] $ \ \textbf{DAZORD \ P., LICHNEROWICZ \
A. \ and \ MARLE \ C.M. }Structure locale des vari\'{e}t\'{e}s de Jacobi.
J.Math. Pures Appl. 70 (1991) 101-152.

\bigskip

$\left[ \text{\textbf{DENL}}\right] $\textbf{\ \ De \ WIDE \ M., LECOMTE \
P. \ }Existence of star products and formal deformations of Poisson algebra
of arbitrary symplectic manifold. Lett. Math Phys. 7 (1983) 487-496

\bigskip

$\left[ \mathbf{DUJ}\right] $ \ \textbf{DUFOUR \ J.P. \ }Normal forms for
Lie algebroids (preprint).

\bigskip

$\left[ \mathbf{GEM}\right] $\textbf{\ \ GERSTENHABER \ M.} \ Deformation of
Rings and Algebras Ann. of Math. 79 (1964) 59-103.

\bigskip

$\left[ \mathbf{GH}\right] $ \ \textbf{GOLDSCHMIDT \ H. \ }Equations de Lie.
J. Diff. Geom 11(1970) 167-223.

\bigskip

$\left[ \mathbf{HEL}\right] $\textbf{\ \ HELMSTETTER \ J. \ }Radical d'une
alg\`{e}bre sym\'{e}trique \`{a} gauche Ann. Inst. Fourier 29 (1979) 17-35.

\bigskip

$\left[ \mathbf{KON}\right] $\textbf{\ \ KONTSEVICH\ M.} \ Deformation
quantization of Poisson manifolds. Prep. q-Alg/970-9040.

\bigskip

\textbf{[KOS] \ KOSMANN-SCHWARZBACH \ Y. \ }Crochet de Schouten et
cohomologies d'alg\`{e}bre de Lie CRAS Paris 312 (1991) R3-126.

\bigskip

\textbf{[KJL}$_{\text{\textbf{1}}}$\textbf{] \ KOSZUL \ J-L. \ }Homologie
des complexes de formes diff\'{e}rentielles d'ordre sup\'{e}rieur Ann.
Scient. Ec. Norm. Sup. 7 (1974) 149-154.

\bigskip

\textbf{[KJL}$_{\text{\textbf{2}}}$\textbf{] \ KOSZUL \ J-L. \ }%
D\'{e}formation des vari\'{e}t\'{e}s localement plates. Ann. Inst. Fourier,
18 (1968) 103-114.

\bigskip

\textbf{[KJL}$_{\text{\textbf{3}}}$\textbf{] \ KOSZUL \ J-L. \ }Crochet de
Schouten-Nijenhuis et cohomologie in Elie Cartan et les Math\'{e}matiques
d'aujourd'hui. Ast\'{e}risque (1985) 253-271.

\bigskip

\textbf{[LMA] \ LIBERMANN \ P., MARLE \ C.M. \ }Symplectic manifolds,
dynamical systems and Hamiltonian Mechanics, in Diff. Geometry and
relativity in Honour of A. Lichn\'{e}rowicz, D. Reidel, Dordrecjt (1976).

\bigskip

\textbf{[MB] \ MALGRANGE \ B. \ }Equations de Lie I,II, Jour. of Diff. Geom.
(6) (1972) 503-522 ep (7) (1972) 117-141.

\bigskip

\textbf{[MJ] \ MILNOR \ J. \ }On the fundamental groups of complete affinely
flat manifolds. Adv. in Math. 25 (1977) 178-187.

\bigskip

\textbf{[NGB}$_{\text{1}}$] \ \textbf{NGUIFFO BOYOM M. \ }Alg\`{e}bres
sym\'{e}triques \`{a} gauche et alg\`{e}bres de Lie r\'{e}ductives.
Th\`{e}se (Grenoble) 1968.

\bigskip

\textbf{[NGB}$_{\text{\textbf{2}}}$\textbf{] \ }Homology of KV-algebras and
related topics in Diff. Geom. (To appear).

\bigskip

\textbf{[NGB}$_{\text{\textbf{3}}}$\textbf{] \ }The Homology theory of
Koszul-Vinberg algebras (submited).

\bigskip

\textbf{[NGB}$_{\text{4}}$] \ Structures affines isotropes. Ann della Sc.
Norm. Sup. Pisa IV (1993) 91-131.

\bigskip

\textbf{[NBW] \ NGUIFFO BOYOM \ N. WOLAK \ R.} \ The KV-homology and
transversally affine foliation (in preparation).

\bigskip

\textbf{[NIJ] \ NIJENHUIS \ A. \ }Sur une classe de propri\'{e}t\'{e}s
communes \`{a} quelques types diff\'{e}rents d'alg\`{e}bres. L'Enseignement
Math\'{e}matiques t. XIV (1969). 225-277.

\bigskip

\textbf{[PAM] \ PEREA \ A.M. \ }Flat left invariant connections adapted to
automorphism structure on Lie groups. Jour. Diff. Geom. 16 (1581). 445-474.

\bigskip

\textbf{[SS] \ SINGER \ I.M. \ and \ STERNBERG \ S. \ }The infinite groups
of Lie and Cartan. Jour. d'Analyse Math. Jerusalem 15 (1965). 1-114.

\bigskip

\textbf{[VAI] \ VAISMAN \ I. \ }Lectures on the Geometry of Poisson
manifolds. Progress in Math. 118 Birkha\"{i}ser 1994.

\bigskip

\textbf{[VEY] \ VEY \ J. \ }D\'{e}formation du crochet de Poisson sur une
vari\'{e}t\'{e} symplectique. Comment. Math. Helv. 50 (1975) 421-454.

\bigskip

\textbf{[VEB] \ VINBERG \ E.B. \ }Theory of convex homogeneous cones, Trudy
Moscow Mat. Obshch. 12 (1963) 303-358.

\bigskip

\textbf{[WA] \ WEINSTEIN \ A. \ }The local structure of Poisson manifolds J.
Diff. Geom. 18 (1983), 523-557.

\end{document}